\documentclass [11pt]{amsart}
\usepackage {amsmath, amssymb, amscd}

\newtheorem{theorem}{Theorem}

\newtheorem {definition}{Definition}
\newtheorem {lemma}{Lemma}
\newtheorem {proposition}{Proposition}

\newtheorem {observation}{Observation}

\def\ss {{\mathcal{S}}}
\def\zz {{\mathbb{Z}}}
\def\rr {{\mathbb{R}}}
\def\cc {{\mathbb{C}}}

\def\tt {{\mathbb{T}}}
\def\bb {{\mathcal{B}}}
\def\inte {{\text{int}}}

\def\zsn {{Z^{\sigma}_{\nu}}}
\def\swf {{\text{SWF}}}  
\def\dir {{\not\!\partial}}
\def\spc {{\mathfrak{c}}}
\def\sps {{\hat{\mathfrak{c}}}}

\def\vp {{\varphi}}   
\def\del {{\partial}}
\def\delt {{\frac{\partial}{\partial t}}}

\def\fin {{\hfill \square}}
\def\sob {{L^2_k}} 
\def\sobo {{L^2_{k+1}}}
\def\sobh {{L^2_{k+1/2}}}
\def\sinh {{\text{sinh}}}
\def\sw {{\text{SW}}}
  
\def\ss {{\textbf{S}}}
\def\ttt {{\textbf{T}}}

\def\pr {{\text{pr}}}
\def\ind {{\text{ind}}}  
\def\coker {{\operatorname{coker}}}
\def\im {{\hspace{2pt} \text{Im }}}

\begin{document}

\title [Gluing theorem] {A gluing theorem for the relative Bauer-Furuta
invariants} \author [Ciprian Manolescu]{Ciprian Manolescu} 

\begin {abstract} In a previous paper we have constructed an invariant of
four-dimensional manifolds with boundary in the form of an element in
the stable homotopy group of the Seiberg-Witten Floer spectrum of the 
boundary. Here we prove that when one glues two four-manifolds along their 
boundaries, the Bauer-Furuta invariant of the resulting manifold is 
obtained by applying a natural pairing to the invariants of the pieces.
As an application, we show that the connected sum of three copies of the 
K3 surface contains no exotic nuclei. In the process we also compute the 
Floer spectrum for several Seifert fibrations.
\end {abstract}

\address {Department of Mathematics, Columbia University\\ 
2990 Broadway, New York, NY 10027}
\email {cm@math.columbia.edu}

\maketitle

\section {Introduction}

In \cite{BF}, Bauer and Furuta have defined an invariant $\Psi$ of closed 
four-manifolds which takes values in an equivariant stable cohomotopy 
group of spheres. As shown by Bauer in \cite{B}, this invariant
is strictly stronger than the Seiberg-Witten invariant. For example, 
$\Psi$ can be used to distinguish between certain connected sums of 
homotopy K3 surfaces; on the other hand, it is well-known that such 
connected sums have trivial Seiberg-Witten invariants. Essential in 
computing the invariant for connected sums was Bauer's gluing 
theorem, which states that the invariant of a connected sum is equal to 
the smash product of the invariants for each of the two pieces.

As a first step in computing the Bauer-Furuta invariants of other
4-manifolds, one would be interested in determining their behavior with
respect to decompositions along more general 3-manifolds. The purpose of
the present article is to generalize Bauer's gluing theorem in this
direction. We will only discuss gluing along rational homology 3-spheres,
but the method of proof seems suitable for further generalizations. In
particular, we do not require any special properties for the metric on the
3-manifold.

In \cite{M} and \cite{KM}, we have extended the definition of $\Psi$ to
compact four-manifolds with boundary. If $X$ has boundary $Y,$ then $\Psi$
takes the form of a stable, $S^1$-equivariant morphism between a Thom
spectrum $\ttt(X)$ associated to the Dirac operator on $X$ and $\swf(Y),$
which is an invariant of 3-manifolds called the Seiberg-Witten Floer
spectrum.  More precisely, $\ttt(X)$ is the formal desuspension of order
$b_2^+(X)$ of the Thom spectrum corresponding to the virtual index bundle
on $H^1(X;  \rr)/H^1(X; \zz)$ coming from the Dirac operators. Note that
all of our invariants depend on the choice of $spin^c$ structures on the
respective manifolds, as well as on choices of orientations for 
$H^2_+,$ but we omit them from notation for simplicity. 

Let $X_1$ and $X_2$ be two compact, orientable four-manifolds with 
boundaries $Y$ and $-Y,$ respectively, where $b_1(Y)=0.$ We can form the 
manifold $X = X_1 \cup_Y X_2$ and it is easy to see that $\ttt(X) = 
\ttt(X_1) \wedge \ttt(X_2).$ Furthermore, we know from \cite{M} that the 
spectra $\swf(Y)$ and $\swf(-Y)$ are Spanier-Whitehead dual to each other. 
There is a natural duality morphism:
$$\eta: \swf(Y) \wedge \swf(- Y) \to \ss,$$
where $\ss$ is the sphere spectrum.

We will prove the following gluing result:

\begin {theorem} 
\label {gluing}
The Bauer-Furuta invariant of $X,$ as a morphism
$$ \Psi(X) : \ttt(X) \to \ss,$$
is given by the formula:
$$ \Psi(X) = \eta \circ \big(\Psi(X_1) \wedge \Psi(X_2) \bigr).$$
\end {theorem}

The main idea in the proof is to do finite dimensional approximation on 
$X_1, X_2,$ and $X$ at the same time, respecting a fiber product formula 
for the respective Sobolev spaces (Lemma~\ref{sobolev}).

In \cite{M} we have noted that the invariant $\Psi$ can also be
interpreted in terms of cobordisms. If $X$ is a 4-dimensional cobordism
between $Y_1$ and $Y_2,$ then the restriction of $\Psi(X)$ to a fiber over
the Picard torus gives a morphism between $\swf(Y_1)$ and
$\swf(Y_2)$ with a possible change in degree, i.e. a morphism
$$ \mathcal{D}_X : \swf(Y_1) \to \Sigma^{b_2^+(X), -d(X)}\swf(Y_2).$$  
Here $\Sigma^{m, n}$ denotes suspension by $m$ real representations and 
$n$ complex ones (this could mean formal desuspension if $n < 0$), while
$ d(X) = (c^2 - \sigma(X))/8, c$ being the first Chern class of the 
determinant line bundle for the $spin^c$ structure on $X.$

For example, in our case $\Psi(X_1)$ gives rise to a morphism between 
$\ss$ and $\swf (Y)$
and, via the duality map, $\Psi(X_2)$ gives a morphism $\mathcal{D}_{X_2}$
between $\swf (Y)$ and $\ss$ (both up to a change in degree).  We can then
rephrase Theorem~\ref{gluing} by saying that the Bauer-Furuta invariant of 
$X$ is the composition of $\mathcal{D}_{X_1}$ with the relevant suspension 
of $\mathcal{D}_{X_2}.$

More generally, the same method of proof applies to compositions 
of cobordisms where the initial and the final 3-manifolds are nonempty. 
Thus we have the following:

\begin{theorem}
\label{glu2}
Let $X_1$ be a cobordism between $Y_1$ and $Y_2,$ and $X_2$ a cobordism 
between $Y_2$ and $Y_3,$ where the $Y_j$'s are homology 3-spheres. Denote 
by $X$ the composite cobordism between $Y_1$ and $Y_3.$ Then, for $b = 
b_2^+(X_1)$ and $d=d(X_1),$ we have:
$$ \mathcal{D}_{X} = (\Sigma^{b, -d}\mathcal{D}_{X_2}) \circ 
\mathcal{D}_{X_1}.$$
\end {theorem}

In order to be able to use these gluing results, in 
sections~\ref{sec:cmh} and \ref{sec:ex} we develop some 
techniques for computing the Seiberg-Witten Floer spectrum for 
rational homology 3-spheres. This computation can be easily carried out 
for elliptic manifolds using a metric of positive scalar curvature. It 
can also be done for some Brieskorn spheres, using the work of Mrowka, 
Ozsv\'ath, and Yu from \cite{MOY}. We describe explicitely $\swf(-Y)$ 
when $Y=\Sigma(2,3,6n \pm 1).$ 

In the last section we give a few topological applications of 
Theorem~\ref{gluing}. When $Y$ is a lens space, for example, we 
recover some versions of the adjunction inequality for Seiberg-Witten 
basic classes. 

When $Y$ is the Brieskorn sphere $\Sigma(2,3,11),$ a form 
of the gluing theorem for relative Seiberg-Witten invariants was 
proved by Stipsicz and Szab\'o in \cite{SS}. Their methods show, for 
example, that the K3 surface does not contain any embedded copy of 
an exotic nucleus $N(2)_{p,q},$ obtained from the standard 
nucleus $N(2)$ by doing logarithmic transformations of 
multiplicities $p$ and $q$ along the fibers, for $p, q \geq 1$ 
relatively prime, $(p,q) \neq (1,1).$ Using Theorem~\ref{gluing}, we 
can strengthen this result:

\begin {theorem}
\label {exotic}
The connected sum $K3 \# K3 \# K3$ does not split off any 
exotic nucleus $N(2)_{p,q},$ for $p, q \geq 1$ 
relatively prime, $(p,q) \neq (1,1).$
\end {theorem}

\noindent{\bf Acknowledgements.} This paper was part of the author's 
Ph.D. thesis, written at Harvard University. I am very grateful to 
my advisor Peter Kronheimer for his guidance and encouragement, as well as 
to Lars Hesselholt, Tom Mrowka, and Cliff Taubes for several helpful 
conversations. 

\section {The setup}

As in the introduction, we let $X_1$ and $X_2$ be two compact, orientable 
four-manifolds with boundaries $Y$ and $-Y,$ respectively. We assume
$b_1(Y)=0,$ and give $X_1$ and $X_2$ 
metrics which are cylindrical near their boundaries. We also endow each 
manifold with a $spin^c$ structure and a base $spin^c$ connection in a 
compatible way, i.e. so that the restrictions of the 4-dimensional objects 
to their boundaries are the 3-dimensional objects. We choose the base 
connection on $Y$ to be flat. 

To distinguish between objects on different manifolds, we will
conventionally use the subscript $j$ the ones on $X_j (j=1,2),$ while 
leaving the three-dimensional ones unmarked. For example, 
we will denote the Dirac operator associated to the 
respective base connection by $\dir_1$ on $X_1,$ $\dir_2$ on $X_2,$ 
and by $\dir$ on $Y.$ The base connections themselves will be 
$A_1, A_2,$ and $A,$ respectively.

\subsection {The Seiberg-Witten Floer spectrum.}
We recall from \cite{M} the definition of the invariant $\swf(Y).$ Let $V$ 
be the space of pairs $x= (a, \phi)$ consisting of a 1-form $a$
and a spinor $\phi$ on $Y,$ with $d^* a=0.$ Coulomb projection is
denoted $\Pi: i\Omega^1(Y) \oplus \Gamma(W) \to V,$ and its linearization
by $\Pi'.$ For $\mu \gg 0
\gg \lambda,$ we consider the finite dimensional approximation 
$V^{\mu}_{\lambda}$ to $V,$ made of eigenspaces of the operator
$(*d, \dir)$ with eigenvalues between $\lambda$ and $\mu.$ We denote by 
$p^{\mu}_{\lambda}$ the orthogonal projection to $V^{\mu}_{\lambda}.$
Similarly, $V^\mu$ (resp. $V_{\lambda}$) will be the direct sum of all 
eigenspaces with eigenvalues at most $\mu$ (resp. at least $\lambda$), and 
by $p^{\mu}$ and $p_{\lambda}$ the corresponding projections. 

The Seiberg-Witten map on $V$ can be written as a sum of a linear Fredholm 
part $l$ and a compact part $c.$ Both $l$ and $c$ are thought of as 
mapping the Sobolev $L^2_{k+1/2}$ completion of $V$ to the $L^2_{k-1/2}$ 
completion of $V,$ where $k > 3$ is an integer. (In the original 
formulation, $k$ was taken to be a half-integer, but everything in 
\cite{M} works for any real $k > 3.$ In Lemma~\ref{sobolev} below we will 
need $k$ to be an integer.)

We can consider the Seiberg-Witten flow on $L^2_{k+1/2}(V^{\mu}_{\lambda})$:
$$ \delt x(t) = -(l + p^{\mu}_{\lambda}c) x(t).$$
We write $x(t) = \varphi_t(x(0))$ for $t \in \rr.$

The following compactness result is proved in \cite{M}:
\begin{proposition}
\label {main}
Given any $R \gg 0,$ and for any $\mu, -\lambda$ sufficiently large 
compared to $R,$ if a flow trajectory $x: 
\rr\to L^2_{k+1/2}(V_{\lambda}^{\mu})$ satisfies $x(t) \in 
\overline{B(2R)}$ for all $t,$ then in fact $x(t) \in B(R)$ for all $t.$
\end {proposition}

This says that the set $S$ of flow trajectories contained in 
$B=\overline{B(2R)}$ is an isolated invariant set, i.e. 
$$ S = \{x \in B: \varphi_t(x) \in B \text{ for all } t \in \rr\} \subset 
\inte B.$$

Given an isolated invariant set for a flow, one can associate to it an 
invariant called the Conley index. It takes the form of a pointed space 
$N/L,$ where $(N, L)$ is a pair of compact subspaces $L \subset N \subset 
B$ with certain properties. There are many suitable index pairs $(N, L),$
but the quotient $N/L$ is well-defined up to canonical homotopy 
equivalence. The spectrum $\swf(Y)$ is then basically the suspension 
spectrum associated to $N/L,$ shifted in dimension by desuspending by 
some finite dimensional subspace. For more details, we refer to \cite{M}.

\subsection {The relative Bauer-Furuta invariants.}
We continue by summarizing the definition of the morphisms $\Psi,$ again 
following \cite{M}.

Let us first introduce some notation which will be useful later on. For 
$j=1,2,$ given a spinor/1-form pair $x_j$ on $X_j,$ we denote by $\phi(x_j) 
\in
\Gamma(W)$ the restriction of the spinor component to the boundary. Then,
we let $\omega(x_j),$ $\beta(x_j)$ and $\gamma(x_j)$ be the $i\ker d^* \in
i\Omega^1(Y),$ $i\ker d \in i\Omega^1(Y)$ and $i\Omega^0(Y) \cdot dt$
components of the 1-form on the boundary, respectively. Note that if the
1-form component of $x_j$ is in the kernel of $d^*,$ then automatically
$\gamma(x_j) \in i\im d^* \subset i\Omega^0(Y).$ We denote collectively
$\alpha = (\omega, \phi)$ and $$i^*= (\alpha, \beta, \gamma).$$

In this notation, for example, Coulomb projection is
$$\Pi i^*(x_j) = (\omega(x_j), e^{d^{-1}\beta(x_j)}\phi(x_j)),$$
with linearization $\alpha = \Pi' i^*.$

Let $\Omega_g^1(X_1)$ be the space of 1-forms $a_1$ on $X_1$ in Coulomb 
gauge: $a_1 \in \ker d^*, a_1|_{\del X_1} (\nu)=0,$ where $\nu$ 
is the normal vector to the boundary. 

The morphism $\mathcal{D}_{X_1}$ arises as the finite dimensional 
approximation of the Seiberg-Witten map:
$$ SW: i\Omega^1_g(X_1) \oplus \Gamma(W^+_1) \to i\Omega^2_+(X_1)
\oplus \Gamma(W^-_1) \oplus V^{\mu}$$ 
$$(a_1, \phi_1) \to  (F^+_{A_1 + a_1} -
\sigma (\phi_1, 
\phi_1), \dir_1 \phi_1 +  a_1 \cdot \phi_1, p^{\mu} \Pi 
i^*(a_1, \phi_1))$$

We will sometimes denote the map to the first two factors as $sw,$ so that 
$SW = (sw, p^{\mu}\Pi i^*).$

The finite dimensional approximation goes as follows. We first extend $SW$
to a map, still called $SW,$ between the $L^2_{k+1}$ completion of the
domain and the $L^2_k$ and $L^2_{k+1/2}$ completions of $i\Omega^2_+(X_1)
\oplus \Gamma(W_1^-)$ and $V^{\mu},$ respectively. We decompose $SW$ as 
$L+C,$ where $L$ is its linearization and
$C$ is compact. Then we pick $U_1$ a finite dimensional subspace of
$\sob (i\Omega^2_+(X_1) \oplus \Gamma(W_1^-)),$ let $U_1' =
L^{-1}(U_1),$ and $$ SW_{U_1} = \pr_{U_1 \times
V^{\mu}_{\lambda}} \circ SW: U_1' \to U_1 \times
 V^{\mu}_{\lambda}.$$
We denote by $B(U_1, \epsilon)$ and $S(U_1,
\epsilon)$ the closed ball and the sphere of radius $\epsilon$ in $U_1$
(in the $\sob$ norm), respectively. Set 
$$\mathcal{M}_1 =SW_{U_1}^{-1}(B(U_1,\epsilon) \times  
V^{\mu}_{\lambda})$$
and let $M_1'$ and $K_1'$ be the intersections of $\mathcal{M}_1$
with $B(U_1', R)$ and $S(U_1', R)$ respectively, for some $R > 0$ 
(in the $\sobo$ norm). We denote by $M_1, K_1$ the images of 
$M_1'$ and $K_1'$ under the composition of $SW_{U_1}$ with projection 
to the factor $V^{\mu}_{\lambda}.$  One can find an
index pair $(N_1,L_1)$ which represents the Conley index for 
$V^{\mu}_{\lambda}$ in the form $N_1/L_1$ such that $M_1 \subset N_1$ and 
$K_1 \subset L_1.$ 

Now we have a map:
$$ SW_{U_1}: B(U_1', R)/S(U_1', R) \to \Bigl(B(U_1, 
\epsilon)/S(U_1, \epsilon)\Bigr) \wedge (N_1/L_1).$$

We will often say that we are interested in properties of this map for 
all $U_1$ {\it sufficiently large.} What this means is that out of 
every nested sequence $(U_1)_n \subset \sob (i\Omega^2_+(X_1) \oplus 
\Gamma(W_1^-))$ such that $\pr_{(U_1)_n} \to 1$ pointwise  
as $n \to \infty,$ these are properties of $(U_1)_n$ for all $n \gg 0.$

In a similar vein, we can choose $U_2$ a finite dimensional 
approximation on $X_2$ and we obtain a map:
$$ SW_{U_2}: B(U_2', R)/S(U_2', R) \to \Bigl(B(U_2,
\epsilon)/S(U_2, \epsilon)\Bigr) \wedge (N_2/L_2).$$

Finally, finite dimensional approximation on $X$ using subspaces $U$ gives 
a map:
$$ SW_U : B(U', R)/S(U', R) \to B(U, \epsilon)/S(U, \epsilon).$$

\subsection {Equivariant maps between spheres.}

The map $SW_U$ can be thought of as an equivariant map between two spheres 
with semifree $S^1$ actions. In general, given such an equivariant map $f: 
E'^+ \to E^+,$ we can choose an identification of $E$ with the standard  
$\cc^{n} \oplus \rr^{m},$ for some $n, m \geq 0.$ Up to homotopy, there 
is a $\zz/2$ ambiguity in this choice, because we can compose with an 
orientation reversing isomorphism. Similarly, we can identify $E'$ with 
some $\cc^{n'} \oplus \rr^{m'}.$ 

If we are given an orientation on $E' \oplus E,$ we choose the isomorphisms 
above so that their direct sum is an orientation preserving map
$$ (E' \oplus E) \to (\cc^{n+n'} \oplus \rr^{m+m'}). $$

In this way $f$ defines a canonical element in the equivariant stable
homotopy groups of spheres. Indeed, due to the additive structure of the
equivariant stable category, if we compose with an orientation reversing
map both on the image and on the target, we get back the same element.
Note that this is not true unstably, as shown by the example of the Hopf 
map $S^3 \to S^2,$ which does not depend on the orientation of the target. 

In our case, the orientation on $U' \oplus U$ is fixed once we specify an 
orientation for $H^2_+(X).$ Similarly, if we fix an orientation for 
$H^2_+(X_1),$ the maps $SW_{U_1}$ give a 
canonical element in the stable homotopy groups of $\swf(Y).$ 

\subsection {From fibers to bundles.}

In \cite{M} we have shown that as we increase $U_1,$ the maps $SW_{U_1}$
give rise to the same element in the stable homotopy groups of $\swf(Y).$ 
These maps depend on the base connection $A_1,$ which we can allow to vary 
by adding harmonic forms to it. The end result is a collection of maps 
parametrized by the Picard torus of $X_1,$ which gives (after 
stabilization) the morphism $\Psi(X_1) : \ttt(X_1) \to \swf(Y).$ Similarly 
we can obtain $\Psi(X_2)$ from a collection of maps parametrized by the 
Picard torus of $X_2.$

In order to prove Theorem~\ref{gluing}, we will focus on showing the 
corresponding statement fiberwise. In other words, we fix $A_1, A_2$ such 
that they both restrict to $A$ on the cylindrical ends, we glue 
them together to give a base connection $A_X$ on $X,$ and then we 
will show:

\begin{theorem}
\label{new}
There exist $\mu, -\lambda, U, U_1, U_2$ sufficiently large and 
$\epsilon >0$ sufficiently small such
that the maps $\eta \circ (SW_{U_1} \wedge SW_{U_2})$ and $SW_U$
are stably homotopic.
\end {theorem}

By ``stably homotopic'' we mean that the two maps induce the same element 
in the equivariant stable homotopy groups of spheres. 

The method of proof can be applied to a bundle of maps over the Picard 
torus rather than to single maps, yielding a proof of 
Theorem~\ref{gluing}. 

\subsection {The duality map.}
\label{sec:duality}
In order to prove Theorem~\ref{new}, we need an explicit description of 
the duality map $\eta.$

The origin of $\eta$ is the duality theorem between the Conley indices for 
the forward and reverse flows (\cite{Cor}). In our situation, the 
flow on $V^{\mu}_{\lambda}$ coming from the Seiberg-Witten flow for $Y$ is 
the reverse of the one coming from the Seiberg-Witten flow for $-Y.$  

According to \cite{McC} and \cite{Cor}, the duality between the Conley 
indices $N_1/L_1$ and $N_2/L_2$ can be represented as 
follows. One can choose the index pairs so that $N_1 = N_2 =N$ 
is a manifold with boundary whose interior is an open subset of 
$V^{\mu}_{\lambda}$ and $\del N = L_1 \cup  L_2,$ where $L_j$ 
are manifolds with boundary $\del L_1 = \del L_2 = 
L_1\cap L_2.$ The duality is represented by a map:
$$ \eta : N/L_1 \wedge N/ L_2 \to B(\epsilon)/S(\epsilon),$$ 
where $B(\epsilon)$ and $S(\epsilon)$ are the ball and the sphere of 
radius $\epsilon$ in $V^{\mu}_{\lambda}.$

Here is one way of constructing $\eta.$ Pick $\tilde N$ any compact subset 
contained in the interior of $N.$ Take a small tubular neighborhood 
$[0, \delta] \times \del N$ of $\del N$ in $N,$ disjoint from 
$\tilde N.$ Then
$$ N_1 = N \setminus \bigl( [0, \delta) \times \inte L_2\bigr)$$
is homotopy equivalent to $N$ via a map $m_1: N \to N_1$ which can be 
chosen to be the identity on $\tilde N$ and such that $m_1(L_1)$ lies in 
the interior of $L_1.$ 

Similarly, 
$$ N_2 = N \setminus \bigl( [0, \delta) \times \inte L_1\bigr)$$
is homotopy equivalent to $N$ via a map $m_2: N \to N_2$ which can be 
chosen to be the identity on $\tilde N$ and such that $m_2(L_2)$ lies in  
the interior of $L_2.$ We can assume that $N_1$ is separated from 
$m_2(L_2)$ by a distance at least $\epsilon,$ and that the same is true 
for $N_2$ and $m_1(L_1).$ We can also assume that $|m_j(x) - x| < 
2\epsilon,$ for any $x \in N$ and $j=1,2.$
 
Then we define:
$$ \eta: N \times N / ((N \times L_2) \cup (L_1 \times N)) \to 
B(\epsilon)/S(\epsilon) $$

$$ \eta(x, y) = \begin{cases}
m_1(x) - m_2(y) & \text{if $|m_1(x)  - m_2(y)| < \epsilon$;}\\
* & \text{otherwise.}
\end{cases}$$

Note that when $x, y \in \tilde N$ and $|x-y| < \epsilon,$ then $\eta(x, 
y) = x-y.$ 

\medskip
Now we have explicit descriptions of both

$$\eta \circ (SW_{U_1} \wedge SW_{U_2}) : \Bigl( B(U_1',
R)/S(U_1', R)\Bigr) \wedge \Bigl( B(U_2', R)/S( U_2', R)  
\Bigr) \to $$ 
$$ B(U_1,\epsilon)/S (U_1, \epsilon)\wedge
B(U_2, \epsilon)/S(U_2, \epsilon) \wedge
\Bigl(B(V^{\mu}_{\lambda}, \epsilon)/S(V^{\mu}_{\lambda}, 
\epsilon)\Bigr)$$
and
$$ SW_U : B(U', R)/S(U', R) \to B(U, \epsilon)/S(U, \epsilon).$$

Note that we have the freedom to choose $U, U_1, U_2$ as we want, 
provided they are sufficiently large.

\section {Compactness of the glued up moduli space}

The definition of the map $\eta$ in section \ref{sec:duality} 
involves the 
homotopy equivalences $m_1$ and $m_2,$ which are somewhat inconvenient to 
deal with. Nevertheless, we know that if $x$ and $y$ lie in a fixed 
compact subset $\tilde K$ in the interior of $N,$ then we have a simple 
description:

$$ \eta(x, y) = \begin{cases}
x - y & \text{if $|x  - y| < \epsilon$;}\\   
* & \text{otherwise.}
\end{cases}$$

It turns out that this description is sufficient when working with the map 
$\eta \circ (SW_{U_1} \wedge SW_{U_2}).$ Indeed, the pairs $(x,y)$ to 
which we apply the map $\eta$ in this case are boundary values of 
connection-spinors pairs $(a_j, \phi_j)$ on $X_j$ such that $sw(a_j, 
\phi_j)$ are small. Furthermore, $\eta(x,y) = *$ unless $|m_1(x) - m_2(y)| 
< \epsilon.$ Since $|m_1(x) - x| < \epsilon$ and $|m_2(y) -y| < \epsilon,$ 
this condition implies $|x - y| < 3\epsilon.$ 

Therefore, we have approximate solutions of the Seiberg-Witten equations 
on both sides with approximate boundary values $x$ and $y$ close to each 
other. The idea is that by gluing together the two solutions we obtain an 
approximate Seiberg-Witten solution on $X,$ whose Sobolev norms we can 
control. In particular, we would also be able to control the size of $x$ 
and $y.$ If we prove that they lie in a compact set in 
$V^{\mu}_{\lambda},$ then by choosing the radius $R$ in the definition of 
$SW_{U_j}$ large enough, we could assume that this compact set is $\tilde 
K \subset N.$ 

The result that we need is Proposition~\ref{glu} below. We denote by 
$\Pi'$ the linearization of the Coulomb projection, i.e. the orthogonal 
projection from $i\Omega^1(Y) \oplus \Gamma(W)$ onto $V.$ 

\begin {proposition}
\label {glu}
There exists a constant $C > 0$ such that, given any $R \gg C,$ for 
every $U_1, U_2, \mu,$ $ -\lambda$ sufficiently large 
and for every $\epsilon$ sufficiently 
small, whenever $x_j \in U_j' (j=1,2)$ satisfy:
$$ \| x_j \|_{\sobo} < R; \| \pr_{U_j} sw(x_j) \|_{\sob} < \epsilon; $$ 
$$ \|p^{\mu}_{\lambda} \Pi i^*(x_1) - p^{\mu}_{\lambda} \Pi 
i^*(x_2)\|_{\sobh} < 3\epsilon,$$
we have the bounds:
$$ \| p^{\mu}_{\lambda} \Pi i^*(x_j) \|_{\sobh} \leq C.$$
\end {proposition}

Note that the requirement $x_j \in U_j'$ automatically implies that:
$$ p^{\lambda} \Pi' i^*(x_1) = p_{\mu} \Pi' i^*(x_2) =0. $$

We start by proving the following:
\begin {lemma}
\label {paris}
There is a constant $C > 0$ such that for any smooth $x_1 \in 
i\Omega^1_g(X_1) \oplus \Gamma(W_1^+)$ we have:
$$ \| p_0 \Pi' i^* (x_1) \|_{\sobh} \leq C \cdot \bigl( \|(d^+ \oplus 
\dir_1)x_1\|_{\sob} + \|x_1\|_{L^2}\bigr).$$
\end {lemma}

\noindent $\textbf {Proof.}$ Let us first prove the statement in the 
particular case when $(d^+ \oplus \dir_1) x_1 =0.$ Consider the 
spinorial part $\phi_1$ of $x_1.$ On the cylindrical end $[-1,0] 
\times Y$ near the boundary, we can write $\phi_1$ as a function $\phi: 
[-1,0] \to \Gamma(W).$ If $\phi_{\lambda}$ are the eigenvectors of 
$\dir$ on $\Gamma(W)$ corresponding to eigenvalue $\lambda,$ with 
$\|\phi_{\lambda}\|_{L^2}=1,$ we have an orthogonal decomposition:
$$ \phi(t) = \sum_{\lambda} c_{\lambda}(t) \phi_{\lambda}.$$

The equation $\dir_1 (\phi_1) =0$ implies $c'_{\lambda}(t) + \lambda 
c_{\lambda}(t) =0,$ so $c_{\lambda}(t) = e^{-\lambda t} c_{\lambda}(0).$
Note that:
$$ \| p_0 \Pi' i^* (\phi_1) \|^2_{\sobh} \sim \sum_{\lambda > 0} 
\lambda^{2k+1} c_{\lambda}(0)^2,$$
while
$$ \| \phi_1\|_{L^2}^2 \sim \sum_{\lambda} \frac{e^{2\lambda} 
-1}{2\lambda}  
c_{\lambda}(0)^2,$$
so clearly the first expression is controlled linearly by the second.

A similar discussion applies to $a_1,$ the differential form part of 
$x_1.$ In this case on the cylindrical end we can decompose $a_1$ as $$ 
a_1=\omega(t) + \beta(t) + \gamma(t)dt, \omega(t) \in \ker d^*, \beta(t)  
\in \ker d, \gamma(t) \in i\Omega^0(Y).$$ We are only interested in 
bounding the positive part of $\omega(0),$ and the operator $d^+$ acts on 
the $\omega(t)$ part as $(\del/\del t)\wedge + *d,$ so indeed we can 
proceed just as we did in the spinorial case.

\medskip
Now consider any smooth $x_1 \in i\Omega^1_g(X_1) \oplus 
\Gamma(W_1^+),$ and let $y_1 = (d^+ \oplus \dir_1) x_1.$

There is an $a \ll 0$ so that the Fredholm map:
$$ (d^+ \oplus \dir_1, p^a \Pi' i^* ): \sobo(i\Omega^1_g(X_1) \oplus 
\Gamma(W_1^+)) \to \sob(i\Omega^2_+(X_1) \oplus \Gamma(W_1^-)) $$
\begin {flushright}
$\oplus \sobh(V^a)$
\end {flushright}
is surjective. 

Thus, we can find $x_1'$ with:
$$ (d^+ \oplus \dir_1) x_1' = y_1, p^a \Pi' i^* (x_1') =0$$
such that we have a bound:
$$ \| x_1'\|_{\sobo} < C' \cdot \|y\|_{\sob}.$$

Since $x_1 - x_1'$ satisfies $(d^+ \oplus \dir_1)(x_1 - x_1') =0,$ we 
have:
$$ \|p_0 \Pi' i^*( x_1 - x_1') \|_{\sobh} \leq C \|x_1 - x_1'\|_{L^2}.$$

Both the $\sobh$ norm of $p_0 \Pi' i^* x_1'$ and the $L^2$ norm of $x_1'$ 
are controlled linearly by the $\sobo$ norm of $x_1',$ and hence by the 
$\sob$ norm of $y.$ We deduce that (for a new constant $C$):
$$ \| p_0 \Pi' i^* (x_1) \|_{\sobh} \leq C \cdot \bigl( \|y\|_{\sob} + 
\|x_1\|_{L^2}\bigr),$$  
as desired.$\fin$

\medskip
\noindent $\textbf {Proof of Proposition~\ref{glu}.}$ Fix $R \gg 0.$ 
Assume the statement is not 
true, so that there exist $\mu_n, -\lambda_n \to \infty, \epsilon_n \to 
0, U_{j,n}$ with $\pr_{U_{j,n}} \to 1,$ and $x_{j,n}$ satisfying the 
conditions above but with 
$$ \| p^{\mu_n}_{\lambda_n} \Pi i^*(x_{j,n}) \|_{\sobh} > C.$$

By Rellich's lemma, we can pick a subsequence of the $x_{j,n},$ still 
denoted as such for simplicity, so that $x_{j,n}$ converges to some $x_j$ 
in the $\sob$ norm. Hence $sw(x_{j,n}) \to sw(x_j)$ in $L^2_{k-1}.$ 

Also, we know that $\pr_{U_{j,n}} sw(x_{j,n}) \to 0$ in $\sob,$ and
$\pr_{U_{j,n}} \to 1,$ so we can deduce that $sw(x_j) = 0.$ Furthermore,
$p^{\mu_n}_{\lambda_n} \to 1$ as well, so $\Pi i^*(x_1) = \Pi i^*(x_2).$
It follows that after adding a closed 1-form to them (i.e. changing the
gauge), we can glue $x_1$ and $x_2$ together to form a solution of the
Seiberg-Witten equations on the closed manifold $X.$ The moduli space of
such solutions is compact, so we get a $\sobo$ bound on their size which
is independent of $R,$ and this persists after changing the gauge back.  
Since we know that $x_{j,n} \to x_j$ in $\sob,$ we
obtain an $\sob$ bound on $x_{j,n}$ and an $L^2_{k-1/2}$ bound on the
boundary values $p^{\mu_n}_{\lambda_n} \Pi i^*(x_{j,n}).$

To obtain the stronger $\sobh$ bound which we need, we proceed as follows.
By the compactness of $sw - (d^+ \oplus \dir_j),$ the $\sob$ bound on
$x_{j,n}$ and the fact that $sw(x_{j,n}) \to 0$ in $\sob,$ we get an 
$\sob$
bound on $(d^+ \oplus \dir_j)x_{j,n}$ which is independent of $R.$ From
Lemma~\ref{paris} it follows that we have such a bound on $\| p_0 \Pi' i^*
(x_{1,n}) \|_{\sobh}.$ If we apply Lemma~\ref{paris} to $X_2$ instead of 
$X_1,$ we get a bound of the same type on $\| p^0 \Pi' i^* (x_{2,n}) 
\|_{\sobh}.$

Now, let us look at the map 
$$\Pi - \Pi'= (0, (e^{d^{-1}\beta} -1) \phi).$$ 
Changing everything on $X_1$ and $X_2$ by the same gauge, we can 
assume $\beta(x_1) = \beta(x_2)=0.$  We know that the $\phi(x_{j,n})$
are bounded in $L^2_{k+1/2}$ norm, while $\beta(x_{j,n})$ converges to $0$ 
weakly in $L^2_{k+1/2},$ implying $d^{-1} \beta(x_{j,n}) \to 0$ strongly in 
$L^2_{k+1/2}.$ From this we get that $(\Pi-\Pi') i^* (x_{j,n} - x_j)  \to 0$ 
in $\sobh,$ which gives a bound  
$$\|(\Pi-\Pi') i^* (x_{j,n})\|_{\sobh} < C',$$ 
with $C'$ independent of $R.$ 

Combining this with the hypothesis, we get:
$$ \|p^{\mu_n}_{\lambda_n} \Pi' i^*(x_{1,n}) - p^{\mu_n}_{\lambda_n} \Pi'
i^*(x_{2,n})\|_{\sobh} < 3\epsilon_n + C'$$

Given the $\sobh$ bounds on $p_0 \Pi' i^* (x_{1,n})$ and $p^0 \Pi' i^*
(x_{2,n})$ obtained above, as well as the hypothesis
$$ p^{\lambda_n} \Pi' i^*(x_{1,n}) = p_{\mu_n} \Pi' i^*(x_{2,n}) =0, $$
we obtain an $\sobh$ bound on $p^{\mu_n}_{\lambda_n} \Pi' i^*(x_{i,n})$ 
and hence on $p^{\mu_n}_{\lambda_n} \Pi i^*(x_{i,n})$ too. This 
completes the proof. $\fin$

\section {Composing homotopies}
\label{sec:htpy}

In this section we present the proof of Theorem~\ref{new}. We need to show
that $\eta \circ (SW_{U_1} \wedge SW_{U_2})$ and $SW_U$ are stably
homotopic.

We will do this by a series of homotopies and identifications of 
domains/targets for different maps. All of the maps that we consider will 
be coming from the following type of construction. We start with a 
continuous map between two separable Hilbert spaces $f: H' \to H$ such 
that 
the zero set $K=f^{-1}(\{0\})$ is compact and such that $f$ decomposes 
into a linear Fredholm part $l$ and a compact part $c.$ (The prototype 
is the Seiberg-Witten map on a closed manifold.) We can 
then conisder finite dimensional approximations $E$ for $H$ and $E' = 
l^{-1}(E)$ for $H'$ and the corresponding maps
$$ f_E = l+ \pr_E c: E' \to E.$$
These might not have compact zero sets. However, if we choose a large 
ball 
$B(H', r)$ containing $f^{-1}(\{0\}),$ it is true that for $E$ 
sufficiently large, the intersection $f_E^{-1}(\{0\}) \cap B(E', r)$ is 
compact. By abuse of language, we will say that $f_E$ has {\bf 
stably compact zero set.} Of course, this definition is meaningful only if 
$f_E$ is part of a collection of finite dimensional approximations 
for a map $f.$ In practice, we will only write down $f_E,$ the way $E$ 
approximates a Hilbert space $H$ being usually self-understood from the 
context. 

For $\epsilon > 0$ small, we can then construct the 
map:
$$ \tilde f_E: B(E', r)/\del B(E', r) \to B(E, \epsilon)/\del B(E, 
\epsilon),$$
by sending $x \in f_E^{-1}(B(E, \epsilon))$ to $f_E(x)$ and everything 
else to the basepoint. Such a map defines an element $e$ in the 
equivariant stable 
homotopy group of spheres once we specify an orientation for $E_1 \oplus 
E_2$ or, equivalently, one for $(\ker l) \oplus (\coker l).$ For 
different $E$'s approximating the same Hilbert space, and for any $R
\gg 0$ and small $\epsilon$, the maps $\tilde f_E$ are 
stably homotopic, in the sense that they define the same element $e.$

All the maps that we want to compare are of the type $\tilde f_E,$ 
but for simplicity we will write down the expressions for $f_E,$ 
and we will say that two $f$ and $g$ are \textbf{stably 
c-homotopic} if $\tilde f$ and $\tilde g$ are stably homotopic.

In this language, given the description of $\eta \circ (SW_{U_1} \wedge 
SW_{U_2})$ in Proposition~\ref{glu}, we
need to show that the following two maps are stably c-homotopic:
\begin {eqnarray}
\label {eins}
U_1' \times U_2' &\to& U_1 \times U_2 \times V^{\mu}_{\lambda}\notag \\
(x_1, x_2) &\to& \bigl(\pr_{U_1}sw(x_1), \pr_{U_2}sw(x_2), 
p^{\mu}_{\lambda} \Pi i^*(x_1 - x_2)\bigr) 
\end {eqnarray}
and
\begin {equation}
\label {zwei}
U' \to U; \hskip7pt x \to \pr_U sw(x). 
\end {equation}

\subsection {Some preliminaries.}
\label{sec:prel}
As we mentioned above, from now on we will only work with maps between
finite dimensional vector spaces approximating Hilbert spaces, and 
such that their zero sets are stably compact. When we say that one 
such map is linear, we imply that the limiting map on Hilbert spaces 
is also linear. 

In comparing such maps we will be using repeatedly the following 
observation.

\begin {observation} Let $f: A \to B$ be continuous and $g: A 
\to C$ linear and surjective such that 
$$ h = (f, g): A \to B \oplus C$$
has stably compact zero set. An orientation on $A \oplus B \oplus C$ 
induces 
one on $\ker g \oplus B$ via the identification $g: (A/\ker g) \to C.$ 

Then 
$h$ and $$f: \ker g \to B$$ are stably c-homotopic. 
\end {observation}

\medskip
\noindent $\textbf {Proof.}$ Let $D$ be the orthogonal complement of $\ker 
g$ in $A.$ Since $g|D: D \to C$ is an isomorphism, we have that $f$ is 
stably c-homotopic to: $$ h': (\ker g) \oplus D \to B \oplus C, 
\hskip7pt h'(x, y) = (f(x), g(y)).$$

On the other hand, $h$ can be written as:
$$ h: (\ker g) \oplus D \to B \oplus C, \hskip7pt h(x, y) = (f(x+y), 
g(y)).$$

We can interpolate from $h$ to $h'$ via the maps $(1-t)h + th', t\in 
[0,1].$ This is possible because the zero sets of all these maps are 
identical, hence stably compact. Indeed, if $x \in \ker g$ and $y \in 
D$, then $g (y) =0$ implies $y=0,$ and from here $(1-t)f(x+y) + tf(x) =0$ 
means simply $f(x) =0.$  $\fin$

\medskip
We will refer to the operation of replacing $h$ with $f|_{\ker g}$ as 
``moving the condition $g=0$ from the map to the domain.''

\medskip
The comparison of $(\ref{eins})$ and $(\ref{zwei})$ will be done in 
several steps.

\subsection {Changing the boundary conditions.}

The conditions $\gamma(x_j) =0$ in the definitions of 
$$U_1' = \{ x_1 \in \sobo(i\Omega^1(X_1) \oplus \Gamma(W_1^+)): (d^+
\oplus \dir_j)x_1 \in U_1; d^*x_1 =0; $$
\begin {flushright}
$p^{\lambda} \alpha(x_1) =0; 
\gamma(x_1) =0 \}$
\end {flushright}
and
$$U_2' = \{ x_2 \in \sobo(i\Omega^1(X_2) \oplus \Gamma(W_2^+)): (d^+
\oplus \dir_2)x_2 \in U_2; d^*x_2=0; $$
\begin {flushright}
$p_{\mu} \alpha(x_2) =0; \gamma(x_2) =0 
\}$
\end {flushright}
make these spaces unsuitable for gluing. 

However, we can construct isomorphisms
$$ \Phi_t : U_1' \times U_2' \to U_t'', t \in [0,1],$$
where
\begin {eqnarray*} U_t'' &=& \{ (x_1, x_2):  x_j \in \sobo(i\Omega^1(X_j) 
\oplus \Gamma(W_j^+)); (d^+ \oplus \dir_j)x_j \in U_j;\\
& & d^* x_j =0; p^{\lambda} \alpha(x_1) =0; 
p_{\mu} \alpha(x_2) =0; \\
& &\gamma(x_1) = \gamma(x_2); (1-t) d(\gamma(x_1) + 
\gamma(x_2)) + t(\beta(x_1) - \beta(x_2)) =0   \}
\end {eqnarray*}
by
$$ \Phi_t (x_1, x_2) = (x_1 + du_1, x_2 + du_2),$$
where $u_j$ are harmonic functions on $X_j$ determined by the boundary 
conditions on $Y$:
$$ u_1' = u_2'; $$
$$(1-t)d(\gamma(x_1) + \gamma(x_2) + u_1' + u_2') + 
t(\beta(x_1) - \beta(x_2) + du_1 - du_2) = 0.$$ 

Here prime denotes the normal derivative at the boundary. Finding the 
harmonic functions above is a mixed Neumann-Dirichlet problem, whose 
solution exists and is unique up to addition of constants. Hence the 
maps 
$\Phi_t$ are well-defined. Their inverses can be easily constructed in a 
similar way, so they are bijections.

Denoting $U'' = U_1'',$ we have that (\ref{eins}) is stably c-homotopic to 
the map:
\begin {eqnarray}
\label {eins1}
U'' &\to& U_1 \times U_2 \times V^{\mu}_{\lambda}\notag \\
(x_1, x_2) &\to& \bigl(\pr_{U_1}sw(x_1), \pr_{U_2}sw(x_2),
p^{\mu}_{\lambda} \Pi i^*(x_1 - x_2)\bigr)
\end {eqnarray}
via the homotopy 
\begin {eqnarray*}
U'_1 \times U'_2 &\to& U_1 \times U_2 \times V^{\mu}_{\lambda}\notag \\
(x_1, x_2) &\to& \bigl(\pr_{U_1}sw(\Phi_t(x_1)), \pr_{U_2}sw(\Phi_t(x_2)),
p^{\mu}_{\lambda} \Pi i^* \Phi_t(x_1 - x_2)\bigr)
\end {eqnarray*}

We need to check that the homotopy goes only through maps whose zero sets
are stably compact. The proof is similar to that of 
Proposition~\ref{glu}. The  
essential point is to check that in the limit $\mu, -\lambda, \sigma, U_1,
U_2 \to \infty,$ the Seiberg-Witten equations on $X_1$ and $X_2$ with 
mixed boundary conditions on $Y$ give a compact moduli space. In our case, 
the boundary conditions are:
\begin {eqnarray}
\label {gauge}
\Pi i^*(x_1)  &=& \Pi i^*(x_2) ; \\ (1-t)(\gamma(x_1) + \gamma(x_2)) &+&
td^*\beta(x_1-x_2)=0;\notag\\ \gamma(x_1) &=& \gamma(x_2) \notag.
\end {eqnarray}  

Let $x_n = (x_{1,n}, x_{2,n})$ be a sequence of such pairs of
solutions. Gauge transform them to $y_n = (y_{1,n}, y_{2,n})$ so that
$y_{1,n} = y_{2,n}$ on $Y.$ Then $y_n$ are actual continuous monopoles on
$X,$ and modulo gauge we can choose a convergence subsequence of them, so
that they are smooth and the convergence is in $C^{\infty}.$ Then, on each
side we can gauge transform them back in a unique way so that they satisfy
($\ref{gauge}$). Since this kind of gauge projection is continuous, we
find that the original subsequence of the $x_n$'s was also convergent.
This implies that the Seiberg-Witten moduli space with boundary conditions 
(\ref{gauge}) is compact.

\subsection {Moving the Coulomb gauge condition from the domain to the 
maps.}

The map (\ref{zwei}) is a finite dimensional 
approximation for 
$$ sw: L^2_{k+1}(i \ker d^* \oplus \Gamma(W^+)) \to L^2_k (i \Omega^2_+(X) 
\oplus \Gamma(W^-)).$$
An alternate approach, which turns out to be more convenient, is to consider 
the map:
$$ s=(sw, d^*): L^2_{k+1}(i \Omega^1(X) \oplus \Gamma(W^+)) \to L^2_k (i 
\Omega^2_+(X) \oplus \Gamma(W^-) \oplus i\Omega^0(X)/\rr).$$

Here $\Omega^0(X)/\rr$ denotes the space of functions which integrate to 
zero on $X.$ Choosing a finite dimensional approximation $T$ for 
$i\Omega^0(X)/\rr,$ we can consider $\tilde U' = (d^+ \oplus \dir \oplus 
d^*)^{-1}(U \times T).$  If $U$ is large enough, then the linear map $d^*: 
\tilde U' \to T$ is surjective and Observation~1 implies that (\ref{zwei}) 
is stably c-homotopic to:
\begin {equation}
\label{zwei1}
\tilde U' \to U \times T; \hskip5pt x \to (\pr_U sw(x), d^*x).
\end {equation}

Starting from here, we can in fact replace (\ref{zwei1}) by other finite 
dimensional approximations to $s=(sw, d^*).$ We can choose any 
sufficiently large subspace $\tilde U \subset L^2_k (i \Omega^2_+(X) \oplus 
\Gamma(W^-) \oplus i\Omega^0(X)/\rr),$ not necessarily 
of the form $U \times T,$ and consider its preimage under the 
linearization $(d^+ \oplus \dir \oplus d^*),$ which we still denote $\tilde 
U'.$ Then we have a map
\begin {equation}
\label{zwei2}
\tilde U' \to \tilde U; \hskip5pt x \to \pr_{\tilde U} s(x),
\end {equation}
stably c-homotopic to the original (\ref{zwei}).

Similarly, we have that (\ref{eins1}) is stably c-homotopic to a map:
\begin {eqnarray}
\label {eins2}
\tilde U'' &\to& \bigl((\tilde U_1 \times \tilde U_2)/\rr\bigr) \times 
V^{\mu}_{\lambda}\notag \\
(x_1, x_2) &\to& \bigl(\pr_{U_1} s(x_1), \pr_{U_2} s(x_2),
p^{\mu}_{\lambda} \Pi i^*(x_1 - x_2)\bigr),
\end {eqnarray}
where $\tilde U_j$ is a finite dimensional approximation for $i
\Omega^2_+(X_j) \oplus \Gamma(W_j^-) \oplus i\Omega^0(X_j),$ the $\rr$ in 
$(\tilde U_1 \times \tilde U_2)/\rr$ stands for 
$$\int_{X_1} d^*x_1 + 
\int_{X_2} d^*x_2 = \int_Y (\gamma(x_1) - \gamma(x_2)) \in \rr,$$
and
\begin {eqnarray*}
\tilde U'' &=& \{ (x_1, x_2):  x_j \in \sobo(i\Omega^1(X_j)  
\oplus \Gamma(W_j^+)); (d^+ \oplus \dir_j \oplus d^*)x_j \in \tilde 
U_j;\\ 
& & p^{\lambda} \alpha(x_1) =0; p_{\mu} \alpha(x_2) =0; 
\gamma(x_1) = \gamma(x_2); \beta(x_1) = \beta(x_2) \}.
\end {eqnarray*}

\subsection {Linearizing the Coulomb projection.}

We change the map (\ref{eins2}) using the homotopy:
\begin {eqnarray}
\label {eins3}
\tilde U'' &\to& \bigl((\tilde U_1 \times \tilde U_2)/\rr\bigr) \times
V^{\mu}_{\lambda}\notag \\
(x_1, x_2) &\to& \bigl(\pr_{U_1}s(x_1),  
\pr_{U_2}s(x_2),
p^{\mu}_{\lambda} ((1-t) \Pi i^* + t\alpha)(x_1 - x_2)\bigr),
\end {eqnarray}
for $t \in [0,1].$ 

To make sure that the zero set is stably compact throughout this 
homotopy, in 
the limit we need to check the compactness of the space of Seiberg-Witten 
equations on $X_1$ and $X_2$ with boundary conditions:
\begin {eqnarray*}
 (1-t)\Pi i^*(x_1) + t \alpha(x_1)  &=& (1-t)\Pi i^*(x_2) + t 
\alpha(x_2);
\\ \beta(x_1)=\beta(x_2) &;& \gamma(x_1) = \gamma(x_2).
\end {eqnarray*}
\noindent Recall that $\alpha(x_j) = (\omega(x_j), \phi(x_j))$ and $\Pi 
i^* (x_j) =$ 
$(\omega(x_j), e^{d^{-1}\beta(x_j)}\phi(x_j)).$ Since $\beta(x_1) =
\beta(x_2),$ the boundary conditions do not actually change throughout 
the homotopy.

\subsection {Moving all boundary conditions from the map to the domain.}

Using Observation~1 again, the map (\ref{eins3}) at $t=1$ is seen to be 
stably c-homotopic to:
\begin {eqnarray}
\label {eins4}
Q' &\to& (\tilde U_1 \times \tilde U_2)/\rr \notag \\
(x_1, x_2) &\to& \bigl(\pr_{U_1}s(x_1),  
\pr_{U_2} s(x_2) \bigr),
\end {eqnarray}
with
\begin {eqnarray*}
Q' &=& \{ (x_1, x_2):  x_j \in \sobo(i\Omega^1(X_j)
\oplus \Gamma(W_j^+)); (d^+ \oplus \dir_j \oplus d^*)x_j \in \tilde
U_j;\\
& & p^{\mu} \alpha(x_1) - p_{\lambda} \alpha(x_2) =0;
\gamma(x_1) = \gamma(x_2); \beta(x_1) = \beta(x_2)) \}.
\end {eqnarray*}

This is true under the hypothesis that the linear map
\begin{equation}
\label {this}
\tilde U'' \to V^{\mu}_{\lambda}; \hskip5pt
(x_1, x_2) \to p^{\mu}_{\lambda} (\alpha(x_1) - \alpha(x_2))
\end {equation}
is surjective.

Set
\begin {eqnarray*}
\tilde Q &=& \{ (x_1, x_2):  x_j \in \sobo(i\Omega^1(X_j)
\oplus \Gamma(W_j^+)); (d^+ \oplus \dir_j \oplus d^*)x_j \in \tilde
U_j;\\ & & \gamma(x_1) = \gamma(x_2); \beta(x_1) = \beta(x_2)) \}.
\end {eqnarray*}
and consider the map 
\begin{equation}
\delta: \tilde Q \to V; \hskip6pt \delta(x_1, x_2) = \alpha(x_1) 
- \alpha(x_2).
\end {equation}
The map $\delta$ is Fredholm. Let us choose $\tilde U_j$ so that $\tilde 
Q$ is sufficiently large for $\delta$ to be surjective.

The following lemma is well-known:
\begin {lemma}
Let $f: H_1 \to H_2$ be a surjective Fredholm map between two Hilbert 
spaces. Let $c_n: H_1 \to H_2$ be a sequence of compact linear maps which 
converge to zero in the operator norm. Then, for $n \gg 0,$ $f_n = f+c_n$ 
is Fredholm and surjective. Furthermore, the orthogonal projections 
$\pr_{\ker f_n}$ converge to $\pr_{\ker f}$ as $n \to \infty.$
\end {lemma}

The maps
\begin{equation}
\label {maps}
\tilde Q \to V; \hskip6pt (x_1, x_2) \to p^{\mu}\alpha(x_1) - 
p_{\lambda}\alpha(x_2)
\end {equation}
are also Fredholm and differ from $\delta$ by compact maps which 
converge to $0.$ Hence for $\mu$ and $-\lambda$ sufficiently large 
compared to $\tilde Q,$ (\ref{maps}) is surjective. Note that this 
automatically implies that (\ref{this}) is surjective. Furthermore, 
the kernel $Q'$ of (\ref{maps}) 
is close to $Q=\ker \delta,$ in the sense that the orthogonal projection 
$ \pr_{Q'} : Q \to V$ is an isomorphism which converges to the 
identity on $Q$ as $\mu \to \infty, \lambda \to 
-\infty.$ 

It follows that (\ref{eins4}) is stably c-homotopic to:
\begin {eqnarray}
\label {drei}  
Q &\to& (\tilde U_1 \times \tilde U_2)/\rr \notag \\
(x_1, x_2) &\to& \bigl( \pr_{U_1} s(x_1),
\pr_{U_2} s(x_2) \bigr),
\end {eqnarray}
with
\begin {eqnarray*}
Q &=& \{ (x_1, x_2):  x_j \in \sobo(i\Omega^1(X_j)
\oplus \Gamma(W_j^+)); (d^+ \oplus \dir_j \oplus d^*)x_j \in \tilde
U_j;\\
& &  \alpha(x_1) - \alpha(x_2) =0;
\gamma(x_1) = \gamma(x_2); \beta(x_1) = \beta(x_2) \}.
\end {eqnarray*}

\subsection {Gluing Sobolev spaces.}

At this point we are left to show that (\ref{drei}) above is stably 
c-homotopic to the map (\ref{zwei2}):
$$ \tilde U' \to \tilde U; \hskip5pt x \to \pr_{\tilde U} s(x).$$

In order to do this, we need to choose $\tilde U_1, \tilde U_2, \tilde U$
carefully so that we can build a suitable identification of some
suspensions of the domains and the targets of the two maps.

We begin by understanding the relationship between the Sobolev spaces
of forms or spinors on $X_1, X_2,$ and $X.$ The relevant result from
analysis is the following:

\begin {lemma} \label{sobolev} Let $X_1, X_2$ be n-dimensional compact
manifolds with common boundary $Y,$ and let $X = X_1 \cup_Y X_2.$ We
assume that $X_1$ and $X_2$ have cylindrical ends near the boundary, and
denote by $t$ the variable in the direction normal to $Y.$ If $E$ is a
vector bundle over $X,$ we denote by $L^2_m(U; E)$ the space of $L^2_m$
sections of $E$ over a subset $U \subset X.$ Then, for every integer $k
\geq 1,$ the space $L^2_k(X;E)$ can be naturally identified as the fiber
product:
$$ L^2_k(X;E) = L^2_k(X_1;E) \times_{\bigl(\prod_{m=0}^{k-1}
L^2_{m+1/2}(Y;E)\bigr)} L^2_k(X_2 ;E)$$ with respect to the maps:
\begin
{eqnarray*} r_j:L^2_k(X_j;E) &\to& \prod_{m=0}^{k-1} L^2_{m+1/2}(Y;E) \\
r_j(u) &=& \Bigl(u|_Y, \frac{\del u}{\del t}\big|_Y, \frac{\del^2 u}{\del
t^2}\big|_Y, \dots, \frac{\del^{k-1} u}{\del t^{k-1}}\big|_Y \Bigr ),
j=1,2. \end{eqnarray*} \end {lemma}

\noindent$\textbf {Proof.}$ It suffices to prove the statement for a neck
surrounding the boundary, so we can assume that instead of $X_1$ and $X_2$
we have $[-1,0] \times Y$ and $[0,1] \times Y,$ respectively, and we
identify $Y$ with $\{0\} \times Y.$

We need to show that if $u_j \in L^2_k(X_j; E) (j=1,2)$ are such that the
restrictions of their $m$th derivatives to $Y$ coincide for $0 \leq m < k,
$ we can combine them to form $u=(u_1, u_2)$ on $X=X_1 \cup X_2$ which is
also in $L^2_k.$

Let us first study the case $k=1.$ Observe that:
$$ L^2_1([a,b] \times Y; E) = L^2_1([a,b], L^2(Y; E)) \cap L^2([a,b],
L^2_1(Y;E)).$$

Since clearly $L^2([-1,1], F) = L^2([-1,0], F) \times L^2([0,1], F)$ for
any $F$ (in particular for $F = L^2_1(Y;E)$), we only have to show
that:
\begin {equation}
\label{unu}
 L^2_1 ([-1,1], F) = L^2_1([-1,0], F) \times_F L^2_1([0,1], F),  
\end {equation}
where $F = L^2(Y; E)$ and the fiber product is taken with respect to the
restrictions to $\{0\}.$ This is well-defined, because for one-dimensional
spaces $L^2_1 \subset C^0.$

In fact, ($\ref{unu}$) is true for any $F.$ Indeed, if $u$ is a continuous
section of $F$ over $[-1,1]$ such that its restrictions to $[-1,0]$ and  
$[0,1]$ are in $L^2_1,$ then the distributional derivative $u'$ is in
$L^2$ when restricted to each of the two halves. Hence $u$ is
differentiable a.e. on each half. Let
$$ v(t) = u(-1) + \int_{-1}^t u'(s)ds. $$

The functions $v$ and $u$ are continuous and their derivatives exist and  
are equal a.e. Furthermore, the distributional derivative $w=(u-v)'$ is  
supported only at $0.$ Note that $w \in L^2_{-1},$ being the derivative of
a continuous function. Hence $w$ is a multiple of the $\delta$ function.
However, $w$ integrates to $0$ on a small interval around $0,$ so in fact
$w=0$ and $u=v.$ Thus $u'$ exists as a function.  Since its restriction
to each half is in $L^2,$ $u'$ itself is in $L^2$ and hence $u \in L^2_1.$
This takes care of the case $k=1.$

For general $k,$ we can apply the $k=1$ statement to $u^{(m)} =   
(u_1^{(m)}, u_2^{(m)})$ for $m < k$ (where the superscript $(m)$ denotes
the $m$th derivative with respect to $t$). We get that $u^{(m)} \in L^2_1$
and, arguing as above, the $u^{(m)}$'s must be the actual distributional
derivatives of $u.$ This implies that $u \in L^2_k.$ $\fin$

\subsection {Finite dimensional approximation on $Y$}

On the cylindrical end, the self-dual forms on $X_j$ are of the form $a 
\wedge dt + *a,$ where $a$ are 1-forms on $Y.$ Thus, the restriction of 
the bundle of self-dual 2-forms to the boundary $Y$ can be identified with 
the bundle of 1-forms on $Y.$ Under the Hodge decomposition, this is 
$i\Omega^1(Y) = i\ker d^* \oplus i\ker d.$ Similarly, the restriction of 
the spinor bundle $W^-$ to $Y$ can be identified with the 3-dimensional 
spinor bundle $W.$

Let us denote 
$$ Z = i\Omega^1(Y) \oplus i\Omega^0(Y) \oplus \Gamma(W).$$
The elements of $Z$ can be written accordingly as $z=(\omega + \beta, 
\gamma, \phi),$ in agreement with the notation introduced in 
Subsection~\ref{sec:prel}. Here $d^*\omega =0$ and $d\beta =0.$  

On a slice of the cylindrical end, we can identify both the source and the 
target of the map $s$ with $Z,$ while the map $s$ itself 
is of the form $\frac{\del}{\del t} + f,$ with $f: Z \to Z$ the 
direct sum of the Seiberg-Witten and $d^*$ maps on $Y.$ The linearization 
of $f$ is:
$$D: Z \to Z, \hskip4pt D(\omega + \beta, \gamma, \phi)= (*d\omega +  
d\gamma, d^*\beta, \dir \phi).$$ 

The fact that this is a self-adjoint linear Fredholm map allows us to do 
finite dimensional approximation on $Z$ using the eigenvalues of $D.$ We 
denote by $Z^{\sigma}_{\nu}$ the direct sum of all eigenspaces of $D$ 
with eigenvalues between $\nu$ and $\sigma,$ where typically $\sigma \gg 
0$ and $\nu \ll 0.$

\subsection {The choice of finite dimensional approximations}
\label{sec:fda}

We choose the subspaces
$$\tilde U_j \subset \sob (i\Omega^2_+(X_j) \oplus  \Gamma(W_j^-) \oplus
i\Omega^0(X_j))
\hskip4pt (j=1,2)$$
and
$$ \tilde U \subset \sob (i\Omega^2_+(X) \oplus \Gamma(W^-) \oplus
i\Omega^0(X)/\rr)$$
so that they are related to each other by a fiber product construction,
which models the one for the Sobolev spaces themselves.

For $j=1,2$ we have restriction maps as in Lemma~\ref{sobolev}:
\begin {eqnarray*}
r_j : \sob (i\Omega^2_+(X_j) \oplus i\Omega^0(X_j) \oplus 
\Gamma(W_j^-))  \to \rr \times \prod_{m=0}^{k-1} L^2_{m+1/2} (Z) \\
x_j \to \Bigl((-1)^{j+1}\int_{X_j} x_j; \hskip4pt x_j|_Y, \frac{\del 
x_j}{\del 
t}\big|_Y, \frac{\del^2
x_j}{\del t^2}\big|_Y, \dots, \frac{\del^{k-1} x_j}{\del t^{k-1}}\big|_Y
\Bigr ).
\end {eqnarray*}
The integral $\int_{X_j} x_j$ refers to integrating the function 
part of $x_j.$

We choose the $\tilde U_j$ so that
\begin {equation}
\label{images}
r_j(\tilde U_j) = \rr \times \prod_{m=0}^{k-1} Z^{\sigma}_{\nu}.
\end {equation}

Once this is true, we let
$$ \tilde U = \tilde U_1 \times_{(\rr \times \prod_{m=0}^{k-1} \zsn)} 
\tilde U_2.$$
Here the fiber product being taken with respect to the maps $r_1$ and 
$r_2.$ 

Given an increasing sequence of $\mu, -\lambda \to \infty,$ we can choose 
the $\tilde U_j$'s as above to form a nested sequence which becomes
sufficiently large. Note that this automatically implies that the
sequence made out of their fiber products $\tilde U$ becomes sufficiently
large on $X.$

\subsection {Rewriting the map (\ref{drei})} Set
$$ Q_j = \{ x_j \in \sobo(i\Omega^1(X_j)
\oplus \Gamma(W_j^+)); (d^+ \oplus \dir_j \oplus d^*)x_j \in \tilde
U_j\}.$$

Then the domain $Q$ of the map $(\ref{drei}):$ 
\begin {eqnarray*} Q &\to&
(\tilde U_1 \times \tilde U_2)/\rr \notag \\ (x_1, x_2) &\to&
\bigl(\pr_{U_1}s(x_1), \pr_{U_2}s(x_2)\bigr)
\end {eqnarray*}
can be expressed as
$$ Q = \{(x_1, x_2) \in Q_1 \times Q_2: i^*(x_1) = i^*(x_2)\}.$$

There is an orthogonal projection
$$ \pr_{\tilde U}: (\tilde U_1 \times \tilde U_2)/\rr \to \tilde U = 
\tilde U_1 \times_{\rr \times \prod_{m=0}^{k-1} \zsn} \tilde U_2.$$

Hence we can replace (\ref{drei}) by
\begin {eqnarray} 
\label {vier}
Q &\to& \tilde U \times \prod_{m=0}^{k-1} \zsn \notag \\ 
(x_1, x_2) &\to&
\bigl(\pr_{\tilde U} (s(x_1), s(x_2)), i^*(s(x_1)^{(m)} - s(x_2)^{(m)}) 
\bigr).
\end {eqnarray}

\subsection {Changing the map (\ref{zwei2})}

In view of Lemma~\ref{sobolev}, the domain $\tilde U'$ 
of the map (\ref{zwei2}):
$$ \tilde U' \to \tilde U; \hskip5pt x \to \pr_{\tilde U} s(x)$$
can be expressed as
$$ \tilde U' = \{(x_1, x_2)\in Q_1 \times Q_2: i^*(x_1^{(m)}) = 
i^*(x_2^{(m)}), m=0,\dots,k \}.$$

\begin {lemma}
If $(x_1, x_2) \in Q,$ then automatically $i^*(x_1^{(m)}) - 
i^*(x_2^{(m)}) \in \zsn$ for $m=0, \dots, k.$ 
\end {lemma}

\noindent $\textbf {Proof.}$ From the definition of $Q,$ we have $(d^+ 
\oplus \dir_j \oplus d^*)x_j \in \tilde U_j.$ In light of (\ref{images}), 
this implies that $i^*\bigl (x_j^{(m+1)} + Dx_j^{(m)}\bigr) \in \zsn$ for 
$m=0, \dots, k.$ Starting from the fact that $i^*(x_1) - i^*(x_2)=0 \in 
\zsn,$ the statement follows easily by induction on $m.$ We are using the 
implication $D(y) \in \zsn \Rightarrow y \in \zsn.$ $\fin$ 

\medskip
Now we can write $\tilde U'$ as the kernel of the linear map
\begin {eqnarray*}
Q &\to& \prod_{m=0}^{k-1} \zsn\\
(x_1, x_2) &\to& i^*(x_1^{(m+1)}) 
- i^*(x_2^{(m+1)}).
\end {eqnarray*}

This map is surjective provided that we chose the $\tilde U_j$'s 
sufficiently large. Using Observation 1 again, (\ref{zwei2}) is stably 
c-homotopic to the map:
\begin {eqnarray}
\label {funf} 
Q &\to& \tilde U \times \prod_{m=0}^{k-1} \zsn \notag \\
(x_1, x_2) &\to&
\bigl(\pr_{\tilde U} (s(x_1), s(x_2)), i^*(x_1^{(m+1)})
- i^*(x_2^{(m+1)}) \bigr).
\end {eqnarray}

\subsection {The final homotopy} We complete the proof of 
Theorem~\ref{new} by exhibiting a homotopy between the maps  
(\ref{vier}) and  (\ref{funf}). We simply choose the linear homotopy:
\begin {eqnarray*}
Q &\to& \tilde U \times \prod_{m=0}^{k-1} \zsn  \\
(x_1, x_2) &\to&
\bigl(\pr_{\tilde U} (s(x_1), s(x_2)), (1-t)(i^*(s(x_1)^{(m)} - 
s(x_2)^{(m)})) \\
& &  + t(i^*(x_1^{(m+1)}) - i^*(x_2^{(m+1)}) \bigr). 
\end {eqnarray*}

Again, we have to show that we have stably compact zero sets throughout 
the
homotopy. In the limit this boils down to proving the compactness of the
space of pairs $(x_1, x_2)$ satisfying:
$$\pr_{L^2_{k}(X)} (s (x_1), s(x_2)) = 0; $$
$$ \Bigl(tx_1^{(m+1)}  + (1-t)s(x_1)^{(m)}\Bigr)|_Y = 
\Bigl(tx_2^{(m+1)}
+ (1-t)s(x_2)^{(m)}\Bigr) |_Y;$$
$$ x_1 |_Y = x_2|_Y. $$
 
Since on the boundary $\del X_j = Y$ we have $s(x_j) = x_j' + f(x_j)$ for
some function $f,$ and $x_1 = x_2$ on $Y,$ by induction on $m$ we have
that $x_1^{(m)} = x_2^{(m)}$ on $Y$ for all $m \leq k,$ which means that
$(x_1, x_2)$ give an element of $L^2_{k+1}(X).$ This implies that
$(s(x_1), s(x_2))$ is in $L^2_k(X),$ so $(x_1, x_2)$ satisfies the
Seiberg-Witten equations on $X.$ Thus, in fact the zero set does not
change throughout the homotopy, and it is the Seiberg-Witten moduli space
on $X.$

\section {Proof of Theorem~\ref{glu2}}

Here we discuss the generalization of Theorem~\ref{new} to the 
gluing of cobordisms. Let $X_1$ be a cobordism between $Y_1$ and $Y_2,$ 
and $X_2$ a cobordism
between $Y_2$ and $Y_3,$ where the $b_1(Y_j)=0$. We call $X$ the composite 
cobordism between $Y_1$ and $Y_2.$ We want to show that the morphism 
$\mathcal{D}_X$ is the composition of $\mathcal{D}_{X_2}$ with the 
relevant suspension of $\mathcal{D}_{X_1}.$

Let us make some simplifications in notation. Let $\bb_j'$ and $\bb_j$ be 
the configuration spaces $i\Omega_g^1 \oplus \Gamma(W^+)$ and $i\Omega^2_+ 
\oplus \Gamma(W^-),$ respectively, for the 4-manifolds $X_j.$ When we talk 
about the analogous spaces for $X,$ we just erase the subscript. Also, 
$V_j$ will stand for the configuration space $i \ker d^* \oplus \Gamma(W)$ 
on the 3-manifold $Y_j,$ while $(V_j)^{\mu}_{\lambda}$ will be the 
corresponding finite dimensional approximations. Finally, $i_j^*$ denotes 
the inclusion of $Y_j$ into one of the 4-manifolds which it bounds.

Then $\mathcal{D}_{X_1}$ is given by a map:
$$ SW_{U_1}: B(U_1', R)/S(U_1', R) \to \Bigl( B(U_1, \epsilon)/ S(U_1, 
\epsilon) \Bigr) \wedge (N_1/L_1) \wedge (N_2/L_2)$$
coming from the finite dimensional 
approximation of:
$$ SW_1 = (sw, p_{\lambda} i_1^*, p^{\mu}i_2^*): \bb_1' \to \bb_1 \oplus 
(V_1)_{\lambda} \oplus (V_2)^{\mu}$$ 
by choosing some linear susbspace $U_1 \subset \bb_1,$ a preimage $U_1'$ 
under the linearized operator, and suitable index pairs $(N_j, L_j)$ for 
the Conley indices of the flows on $(V_j)^{\mu}_{\lambda}, j=1,2.$
Note that this gives an element in a stable homotopy group of  
$\swf(-Y_1) \wedge \swf(Y_2).$ Via the duality map on the first 
factor, in the stable homotopy category 
this is equivalent to giving a morphism $\mathcal{D}_{X_1}$ between 
$\swf(Y_1)$ and the relevant suspension of $\swf(Y_2).$ 

Similarly, $\mathcal{D}_{X_1}$ is a map
$$ SW_{U_2}: B(U_2', R)/S(U_2', R) \to \Bigl( B(U_2, \epsilon)/ S(U_2,
\epsilon) \Bigr) \wedge (N_2/\bar L_2) \wedge (N_3/L_3)$$
coming from the finite dimensional 
approximation of:
$$ SW_2 = (sw, p_{\lambda} i_2^*, p^{\mu}i_3^*): \bb_2' \to \bb_2 
\oplus (V_2)_{\lambda} \oplus (V_3)^{\mu}. $$

Here the index pair $(N_2, \bar L_2)$ is chosen so that $\del N_2 = L_2 
\cup \bar L_2$ and $\del L_2 = \del \bar L_2 = L_2 \cap \bar L_2.$

To get the composition $\mathcal{D}_{X_2} \circ \Sigma^* 
\mathcal{D}_{X_1},$ we smash $SW_{U_1}$ and $SW_{U_2}$ and then apply the 
duality map $\eta_2$ on the $(N_2/ L_2) \wedge (N_2 / \bar L_2)$ factor.  
As in Section 3, we would like to replace this duality map by its linear 
version, i.e. taking the difference in $(V_2)^{\mu}_{\lambda}$ when the 
two elements are within $\epsilon$ distance of each other. To do this, we 
need an analogue of Proposition~\ref{glu}.   

In the case of closed $X$, Proposition~\ref{glu} is based on the 
compactness of the Seiberg-Witten moduli space on $X.$ When $X$ has 
boundary, the best we can hope is that when we add a cylindrical end
to each component of the boundary, the monopoles with finite energy on the 
resulting manifold are bounded in a suitable norm. This corresponds to the 
fact that the restrictions of the monopoles on $X$ to the boundary $Y = 
Y_1 \cup Y_3$ can be used to construct a map to the Conley index of the 
Seiberg-Witten flow $\varphi$ on 
$$V^{\mu}_{\lambda} = (V_1)^{\mu}_{\lambda} 
\oplus (V_3)^{\mu}_{\lambda}.$$ Recall that the Conley index is the 
quotient $N/L = (N_1/L_1) \wedge (N_3/L_3)$ for an index pair $(N, L)$ 
considered with respect to the isolated invariant set $S \subset 
V^{\mu}_{\lambda}.$ 

Let us make the following 
\begin{definition}
A pair of compact sets $(K_1, K_2)$, $K_2 \subset K_1 \subset 
V^{\mu}_{\lambda},$ is called a {\bf pre-index pair} if for $R \gg 0$ 
the following are satisfied:

1. $K_1$ is contained in the closed ball $B=\overline{B(2R)}$ in
$V^{\mu}_{\lambda}$ taken in the $L^2_{k+1/2}$ norm, and $R$ is large
enough to work in Proposition~\ref{main};

2. If $x \in K_1$ satisfies $\varphi_t(x) \in B$ for all $t > 0,$ then 
$\varphi_t(x) \not \in \del B$ for any $t > 0;$

3. For every $x \in K_2$ there exists $t \geq 0$ such that $\varphi_t(x) 
\not 
\in B.$ 
\end {definition}

According to Theorem 4 in \cite{M}, given a pre-index pair $(K_1, K_2),$ 
one can find an actual index pair $(N, L)$ for $S$ with $K_1 \subset N$ 
and $K_2 
\subset L.$

For $R_0, C > 0,$ we let $K_1$ be the set of pairs 
$(p^{\mu}_{\lambda} \Pi i_1^*(x_1), p^{\mu}_{\lambda} \Pi i_3^*(x_2))   
\in V^{\mu}_{\lambda}$ coming from $x_j \in U_j' (j=1,2)$ which satisfy:
$$ \| x_j \|_{\sobo} < R_0; \| \pr_{U_j} sw(x_j) \|_{\sob} < \epsilon; $$  
$$ \|p^{\mu}_{\lambda} \Pi i_2^*(x_1) - p^{\mu}_{\lambda} \Pi
i_2^*(x_2)\|_{\sobh} < 5\epsilon,$$
and let $K_2$ be the subset of $K_2$ of pairs coming from $x_j \in U_j'$ 
which also satisfy:
$$ \| p^{\mu}_{\lambda} \Pi i_2^*(x_j) \|_{\sobh} \geq C$$
for $j=1$ or $2.$

The analogue of Proposition~\ref{glu} is then:
\begin {proposition}
\label {glucobordisms}
There exists a constant $C > 0$ such that, given any $R_0 \gg C,$ for 
every $U_1, U_2, \mu,$ $ -\lambda$ sufficiently large 
and for every $\epsilon$ sufficiently small, $(K_1, K_2)$ is a 
pre-index pair.
\end {proposition}

The proof is similar to that of Proposition~\ref{glu}, so we omit the 
details.

Next, instead of comparing maps with stably compact zero sets via stable 
c-homotopies as we did in Section~4, we need to introduce the following 
notion:

\begin {definition}
Given two finite dimensional vector spaces $E, E'$ with a chosen 
orientation for $E \oplus E',$ a map 
$$f = (g, h): E' \to E \times V^{\mu}_{\lambda}$$
is called {\bf of Conley type} if 
$$ K_1 = h \bigl(g^{-1}(B(E, \epsilon)) \cap B(E', r) \bigr) $$
and
$$ K_2 = h \bigl (g^{-1}(B(E, \epsilon)) \cap \del B(E', r) \bigr)$$
form a pre-index pair, for $\epsilon > 0$ small and any $r \gg 0.$ 
\end {definition}

This is the analogue of maps with compact zero sets. In practice, our 
maps come as finite dimensional approximations to some map between 
Hilbert spaces, and they are only {\bf stably of Conley type}. This 
means that, rather than $(K_1, K_2)$ being a pre-index pair for every 
$r \gg 0,$ what happens is that for every fixed $r \gg 0,$ when $E, 
E'$ are large enough 
approximations, $(K_1, K_2)$ is a pre-index pair for that $r.$

Given a map $f$ that is stably of Conley type, we can find an 
index pair $(N, L)$ with $K_1 \subset N, K_2 \subset L$ and construct an 
element in the stable homotopy groups of the Conley index $N/L:$
$$ \tilde f: B(E', r) / \del B(E',r) \to B(E, \epsilon) \times N/(S(E, 
\epsilon) \times N \cup B(E, \epsilon) \times L). $$

We say that two maps $f_1, f_2$ are {\bf stably Conley c-homotopic} if 
the corresponding maps $\tilde f_1, \tilde f_2$ are stably homotopic.

In this language, we can rephrase Theorem 2 by taking into account the 
result of Proposition~\ref{glucobordisms}:

\begin {proposition}
The maps
\begin {eqnarray*}
U_1' \times U_2' &\to& U_1 \times U_2 \times (V_2)^{\mu}_{\lambda} 
\times V^{\mu}_{\lambda} \\
(x_1, x_2) &\to& \bigl(\pr_{U_1}sw(x_1), \pr_{U_2}sw(x_2),
p^{\mu}_{\lambda} \Pi i_2^*(x_1 - x_2), 
\end{eqnarray*}
\begin {flushright}
$ p^{\mu}_{\lambda} \Pi 
(i_1^*(x_1), i_3^*(x_2)) \bigr)$
\end {flushright}
and
\begin {equation*}
U' \to U \times V^{\mu}_{\lambda}; \hskip7pt x \to \bigl(\pr_U sw(x), 
p^{\mu}_{\lambda} \Pi (i_1^*(x_1), i_3^*(x_2)) \bigr)
\end {equation*}
are stably Conley c-homotopic.
\end {proposition}

For the proof, all the arguments in Section~4 carry over, with 
stable c-homotopies being replaced by stable Conley c-homotopies.

\section {Computations in Morse homotopy}
\label{sec:cmh}

In order to be able to make use of Theorem~1, we should first develop a 
technique for computing the Floer spectra $\swf(Y)$ for various homology 
spheres $Y.$ In this section we outline the general method for 
calculating Conley indices for gradient flows, and in the next one we use 
it to calculate $\swf(Y)$ for some Seifert fibrations.

\subsection {Morse homotopy in finite dimensions.} \label{sec:ms} Let us 
briefly review 
the decomposition of the Conley index for a Morse-Smale downward gradient 
flow $\varphi$ on a finite dimensional manifold $M$. Our exposition is 
inspired 
from \cite{CJS}, \cite{Con}, and \cite{F}.

We assume that $M$ is stably parallelizable. (In fact, in our applications 
$M$ will always be a vector space.) Let $S$ be a compact isolated 
invariant subset of $M$ consisting of
finitely many critical points $x_1, \dots, x_n$ and flow trajectories
between them. Let $E = \{x_1, \dots, x_n\}.$ For each $x 
\in E,$ we can define its stable
and unstable manifolds: $$ W^u(x) = \{ p \in S: \lim_{t \to -\infty}
\varphi_t(p) = x \}, $$
$$ W^s(x) = \{ p \in S: \lim_{t \to \infty} \varphi_t(p) = x \}. $$

The {\bf Morse-Smale condition} says that each critical point $x \in E$ 
is nondegenerate of some Morse index $\mu(x) \in \zz,$ and that $W^u(x)$ 
intersects $W^s(y)$ transversely in a $(\mu(x) - \mu(y))$-dimensional 
manifold $M_{xy},$ for $x, y \in E.$ The manifold $M_{xy}$ has a natural 
compactification $\overline{M}_{xy},$ given by adding broken trajectories 
going through critical points with index between $\mu(y)$ and $\mu(x).$ This 
$\overline{M}_{xy}$ is a manifold with corners. 

Denote 
$$ E_n = \{x \in E: \mu(x) = n \}.$$
Then $E_n$ is an isolated invariant set with an associated Conley index 
$I(E_n).$ This is a wedge of $n$-dimensional spheres, one for each $x \in 
E_n.$

Let $S_n$ be the union of the unstable manifolds $W^u(x)$ for all $x 
\in E$ with $\mu(x) \leq n.$ Each $S_n$ is also an isolated invariant set, 
and has  
a certain Conley index $I(S_n).$ If $m$ is the maximal value of 
all $\mu(x),$ then $S_m = S.$ 

Our goal is to reconstruct the stable homotopy type of $I(S)$ from 
the topology and the framing of the manifolds of flow lines 
between critical points. We do this for each $I(S_n),$ inductively on $n.$

We begin with $S_0 = E_0,$ which is just a finite union of critical 
points. Its Conley index is a corresponding wedge of 0-spheres, i.e. a 
finite CW complex of dimension 0.

Next, $S_0$ is an attractor subset of $S_1,$ in the sense that all points in 
a small neighborhood of $S_0$ in $S_1$ are taken to $S_0$ by the flow, 
asymptotically as $t\to \infty.$ In fact, the only points in $S_1$ not 
taken to 
$S_0$ are those in $E_1.$ $E_1$ is called the repeller set and there is 
a Puppe sequence in homotopy relating the Conley indices:
$$ I(S_0) \to I(S_1) \to I(E_1) \to \Sigma I(S_0) \to \Sigma I(S_1) \to 
\cdots $$

This gives $I(S_1)$ the structure of a 1-dimensional CW complex, coming 
from attaching to $I(S_0)$ one 1-cell for each element of $E_1.$ The 
homotopy type of the suspension $\Sigma I(S_1)$ is then that of the cone 
of the connecting map:
$$ f: I(E_1) \to \Sigma I(S_0).$$  
This is the suspension of the attaching map. It consists of a collection 
of pointed maps between 1-spheres $I(\{x\}) \to 
\Sigma I(\{y\})$ for $x\in E_1, y \in E_0.$ The degree of such a map (with 
appropriate orientations) is the signed count of the gradient flow lines 
between $x$ and $y.$ 

In general, to go from $I(S_{n-1})$ to $I(S_n)$ we use the 
attractor-repeller sequence:
$$ I(S_{n-1}) \to I(S_n) \to I(E_n) \to \Sigma I(S_{n-1}) \to \Sigma 
I(S_n) \to \cdots $$
It follows that $I(S_n)$ has the structure of an $n$-dimensional CW 
complex. The connecting map is the wedge of suspensions of attaching maps:
$$ f_x : S^n \to \Sigma I(S_{n-1}),$$ 
for each $x \in E_n.$ 

The simplest case is when $S_{n-1} = E_j = \{y\}$ for some $j < n-1.$
Then $f_x$ is an element in the $(n-j-1)$th stable homotopy group of
spheres. Franks proved in \cite{F} that, under the Pontrjagin-Thom
construction, this element corresponds to the framed $(n-j-1)$
dimensional manifold of flow lines $M_{xy}/\rr.$ (\cite{F}) The framing 
comes from
trivializing the normal bundle of the contractible manifold $W^s(y)
\subset M,$ and then restricting this framing to $M_{xy},$ viewed as a
submanifold of $W^u(x).$

More generally, Cohen, Jones, and Segal have shown how the stable
homotopy class of $f_x$ can be recovered from the topology of the flow
using a framed topological category. For more details of their
construction, we refer the reader to \cite{CJS}. In our examples, however, we
will only need to consider attaching maps between cells of consecutive
dimensions. In fact, we make the following definition:

\begin {definition}
An isolated invariant set $S$ as above is called {\bf simple} if there is 
some integer $p$ such that $E_n = \emptyset$ unless $p=n$ or $p=n+1.$
\end {definition}

Given a simple isolated invariant set, its Conley index is the cone of a
map between the corresponding wedges of $p$- and $(p+1)$-spheres, and 
this map is determined by the count of gradient flow lines.

\subsection {Morse homotopy in infinite dimensions.}

In practice, to compute the Seiberg-Witten Floer spectrum one starts 
with a description of the critical points and flow lines in infinite 
dimensions. Assuming the Morse-Smale condition, in good cases one can 
compute a stable 
homotopy type using the procedure outlined above. The Floer spectrum, 
however, is defined from the finite dimensional approximations, so we 
need to check that it gives back the same information.

To do this, we consider certain perturbations of the 
Seiberg-Witten flow which affect neither the possibility of finite 
dimensional approximations nor its result, the spectrum $\swf(Y).$ 
The following definition is useful:

\begin {definition}
Let $E$ be a bundle over a finite dimensional manifold. We denote by 
$V=\Gamma(E)$ the space of smooth sections of $E$ and by $L^2_{k}(V)$ its 
$L^2_k$ Sobolev completion. A map $c:V \to V$ is called 
{\bf very compact} if, for every $k > 3,$ $c$ extends to a compact 
map:
$$ c: L^2_{k+1}(V) \to L^2_k(V) $$
and its differentials extend to continuous maps:
\begin {eqnarray*}
 L^2_{k+1}(V) \times L^2_{k-i_1}(V) \times \dots L^2_{k-i_j}(V) 
&\to& L^2_{k-(i_1 + \dots +i_j)}(V),\\
(x; v_1, \dots v_j) &\to& (d^jc)_x(v_1, \dots, v_j).
\end {eqnarray*}
\end {definition}

The prototype of a very compact map is a quadratic one, or the nonlinear
part of the Seiberg-Witten map. The motivation for introducing this
notion is that the properties listed above were the ones used in \cite{M}
to prove Proposition~\ref{main} for approximate Seiberg-Witten flow
trajectories.  More generally, the same method of proof gives the
following:

\begin {proposition}
\label{maine}
Let $E$ be a bundle over a finite dimensional
manifold with $V= \Gamma(E)$, $l: V \to V$ a self-adjoint, linear,
elliptic differential operator of order one, and $c: V \to V$ a very
compact map. We can approximate $V$ by finite dimensional subspaces 
$V^{\mu}_{\lambda}$ using the eigenspaces of $l$ with eigenvalues between 
$\lambda$ and $\mu.$ For $k > 3,$ let $N$ be a bounded, closed subset 
of $L^2_{k+1}(V)$ such that all flow trajectories $x : \rr \to 
L^2_{k+1}(V),$
$$ \delt x(t) = -(l+c)x(t) $$
which lie inside $N$ are in fact contained in $V \cap U,$ for a fixed 
open subset $U \subset N.$  

Then, for every $\mu, -\lambda \gg 0,$ if an approximate trajectory $x : 
\rr \to L^2_{k+1}(V^{\mu}_{\lambda}),$
$$ \delt x(t) = -(l+p^{\mu}_{\lambda}c)x(t), $$
is contained in $N,$ it must be contained in $U.$
\end {proposition}

In the situation described in the proposition, we can define a Conley 
index $I^{\mu}_{\lambda} = I(S^{\mu}_{\lambda})$ for the set 
$S^{\mu}_{\lambda}$ of 
approximate flow trajectories contained in $U \cap V^{\mu}_{\lambda}.$ 
The suspension spectrum 
$$I = \Sigma^{-V^0_{\lambda}} I^{\mu}_{\lambda}$$ 
is then independent of $\lambda$ and $\mu,$ up to canonical equivalence. 
(See \cite{M} for the model proof in the Seiberg-Witten case.) 

On the other hand, we can also look at the original flow on $N \subset 
L^2_{k+1}(V)$ and at the
set $S$ of flow trajectories inside $N.$ Assuming that this is a gradient
flow, we say that its fixed points are nondegenerate if the Hessian of the
respective functional is nondegenerate at those points. We can define the
stable and unstable Hilbert manifolds of critical points as before, and
say that the flow is Morse-Smale if they all intersect transversely in
finite dimensional manifolds. In this situation, we can define a relative
Morse index of critical points as in Floer theory. In principle, we 
should be able to calculate a
stable homotopy type from the topology of spaces of flow lines in $S$ as
in Subsection~\ref{sec:ms}. However, defining suitable framings for these 
spaces is a nontrivial problem. For this reason, we limit our discussion 
to simple flows, when we have only critical points of two successive 
indices. We can then define a stable homotopy tyep as before, by 
taking two wedges of spheres and counting gradient flow lines. We denote 
this stable homotopy type 
by $I(S),$ and call it the {\bf Morse stable homotopy type of $S.$} 
As it stands, it is only defined up to suspension. 

The following proposition shows how to get hold of $I$ when we only have 
information about the flow on $V:$ 

\begin {proposition}
\label {mhom}
If the infinite dimensional flow on $N$ is Morse-Smale and $S$ is simple, 
then $I$ and $I(S)$ represent the same stable homotopy type, up to 
suspension.
\end {proposition}

\noindent {\bf Proof. } The idea is to deform the flow on $N$ to an 
equivalent one on a finite dimensional approximation, where we know that 
the Conley index computes Morse stable homotopy.

Let us warm up by explaining how to do this in a simple case, when $S$ 
consists of just one critical point $y,$ so that $I(S)$ is a sphere. For 
$\mu, -\lambda \gg 0,$ let 
$y^{\mu}_{\lambda}$ be the orthogonal projection of $y$ onto 
$V^{\mu}_{\lambda}.$ Set
$$ v^{\mu}_{\lambda} =  y - y^{\mu}_{\lambda}.$$ 
Note that $v^{\mu}_{\lambda} \to 0$ as $\mu, -\lambda \to \infty.$ For 
some sufficiently large $\mu$ and $\lambda,$ let us consider an open 
convex subset $U'$ of $V$ with $y \in U' \subset U$ and $y^{\mu}_{\lambda} 
\in U' \subset (U - v^{\mu}_{\lambda}),$ the latter being the translation 
of $U$ by $v^{\mu}_{\lambda}.$

If the original flow $\varphi$ was given by $l+c,$ consider now the flow 
$\tilde \varphi:$
$$ \delt x = -(l+c)(x(t) + v^{\mu}_{\lambda}) = -(l+\tilde c)x(t),$$
where the map
$$ \tilde c(x) = c(x+v^{\mu}_{\lambda}) + l(v^{\mu}_{\lambda})$$
is easily seen to be very compact.

We can apply Proposition~\ref{maine} to the flow $\tilde \varphi$ on
$\tilde N = (N- v^{\mu}_{\lambda}).$ In the infinite dimensional space,
this is just a translation by $v^{\mu}_{\lambda}$ of the flow $\varphi.$
However, their finite dimensional approximations are different. What
happens is that now $\tilde y= y^{\mu}_{\lambda}$ is a nondegenerate
critical point of $\tilde \varphi$ which is contained in
$V^{\mu}_{\lambda}.$ It follows that $\tilde y$ remains a nondegenerate
critical point for the flow of $l + p^{\mu'}_{\lambda'}\tilde c,$ for
every $\mu' \geq \mu, \lambda' \leq \lambda.$ Furthermore, for $\mu',
-\lambda' \gg 0,$ it is the only critical point of the approximate flow
in $U'.$ One can see this, for example, by applying
Proposition~\ref{maine} for $\tilde N$ and taking as open subset an
arbitrarily small ball around $\tilde y.$ Therefore, for such $\mu'$ and
$\lambda',$ the Conley index of $\tilde S^{\mu'}_{\lambda'}=\{\tilde x\} 
\subset U'$ is a sphere. Moreover, on $U'$ we can connect $\varphi$ and 
$\tilde \varphi$ by a homotopy of flows $t\varphi + (1-t)\tilde \varphi, 
t\in
[0,1].$ The usual deformation argument in Conley index theory gives that 
$I(S^{\mu'}_{\lambda'})
= I(\tilde S^{\mu'}_{\lambda'})$ is also a sphere, which exactly matches
$I(S).$

In the general case, $S$ is some finite dimensional stratified space. 
We seek to replace the flow $\varphi$ with a flow $\tilde \varphi$ of the 
form:
$$ \delt x = -(l+c)(x(t) + h(x(t))) = -(l+\tilde c)x(t),$$
where $\tilde c$ is a very compact map of the form
$$ \tilde c(x) = c(x + h(x)) + l(h(x)).$$

We would like the flow $\varphi$ in a neighborhood of $S$ to look the same 
as the flow $\tilde \varphi$ in a neighborhood of an isolated invariant 
set $\tilde S$ which is contained in some finite dimensional approximation 
$V^{\mu}_{\lambda}.$ It could happen that for any $\mu, -\lambda \gg 0,$ 
the orthogonal projection of $S$ onto $V^{\mu}_{\lambda}$ is not 
diffeomorphic to $S,$ so we cannot simply set $\tilde S = 
p^{\mu}_{\lambda}(S).$ However, $S$ is finite dimensional, so a simple 
transversality argument shows that, for $\mu$ and $-\lambda$ sufficiently 
large, we can find a linear subspace $W$ with the following properties: 
$V^{\mu}_{\lambda}$ is the image of $V$ under some orthogonal 
transformation $A$ close to the identity, and the orthogonal projection 
$\text{pr}_W: S \to \hat S \subset W$ is a diffeomorphism. Let $\tilde S 
\subset V^{\mu}_{\lambda}$ be the image of $\hat S$ under $A.$ We can 
arrange so that the segments joining each $x \in S$ to $(A \circ 
\text{pr}_W)(x) \in \tilde S$ form an interval bundle $B$ over $S,$ 
smoothly embedded in $L^2_{k+1}(V).$

Take a small tubular neighborhood $T$ of $S$ in $L^2_{k+1}(V),$ in the
form of a disc bundle $\pi: T \to S$ such that $B$ is its subbundle. A
fiber of $\pi$ is a ball in a finite codimension linear subspace of
$L^2_{k+1}(V).$ Note that $\tilde S \subset T$ and $(\pi \circ A \circ
\text{pr}_W) (x) = x$ for all $x \in S.$ We can define a smooth ``bump''
function $\beta: L^2_{k+1}(V) \to [0,1]$ such that $\beta \equiv 1$ in a  
smaller neighborhood $T'$ of $B,$ $T' \subset T,$ and $\beta \equiv 
0$ outside $T.$ Then we define the functions
$$ \gamma : L^2_{k+1}(V) \to R=(S \times [0,1])/(S \times \{0\}) $$
by
$$ \gamma (x) = \begin{cases}
(\pi(x), \beta(x)) & \text{ if } x \in T; \\
* & \text{ otherwise} 
\end {cases}$$
and
$$ \tau: R \to L^2_{k+1}(V), \hskip6pt \tau(x, r) = r \cdot \bigl( x - (A 
\circ \text{pr}_W) (x) \bigr).$$ 

Define $h$ to be the composition 
$$ h= \tau \circ \gamma: L^2_{k+1}(V) \to L^2_{k+1}(V).$$

With this definition of $h$, the very compactness of the map $\tilde c$
is a consequence of the very compactness of $c,$ together with the fact
that $h$ (and its derivatives) factor through maps to the finite
dimensional space $R.$

Furthermore, up to suspension, $I(S)$ for the flow $\varphi$ is the same 
as
the Conley index $I(\tilde S)$ for the flow
$$ (\del/\del t) x(t) = -(l + p^{\mu'}_{\lambda'}\tilde c)x(t)$$ 
approximating $\tilde \varphi$ on $V^{\mu'}_{\lambda'}$ for $\mu' \geq 
\mu, \lambda' \leq \lambda.$ (The 
flows $\vp, \tilde vp$ look the same in small neighborhoods of $S$ and 
$\tilde S,$ respectively.) Note that $\tilde S$ is the invariant part of 
$T' \cap
V^{\mu'}_{\lambda'}$ in the flow generated by $l+p^{\mu'}_{\lambda'}\tilde
c.$ By the continuation argument, its Conley index is the same as that of
the invariant part $S^{\mu'}_{\lambda'}$ of $T' \cap V^{\mu'}_{\lambda'}$
in the flow 
$$(\del/\del t) x(t) = -(l + p^{\mu'}_{\lambda'}c)x(t).$$ 
The latter represents the spectrum $I,$ and this completes the proof.
\hskip6pt $\fin$
\medskip

Three remarks are in place. First, the result above should be valid
without the hypothesis that $S$ is simple, provided that one could define
a good framed topological category in the infinite dimensional situation,
in the spirit of \cite{CJS}. Second, the result can be easily extended to
$G$-equivariant Morse-Bott-Smale flows when $G$ is a compact Lie group,
under the hypothesis that the stable homotopy type $I(S)$ can be computed
from counting gradient flow lines. (We will see such examples in the next
section.) Third, sometimes we might not have enough information about the
flow on $V$ to compute the whole Morse stable homotopy type.  However, we
might have some partial information, e.g. the knowledge of the number of
critical points and the number of flow lines between points of
consecutive indices, which gives Morse-Floer homology. It may also happen
that the original flow itself is not Morse-Smale. Nevertheless, the
principle of Proposition~\ref{mhom} remains true: whatever information we
get out of the infinite-dimensional flow is the same as the information
coming from the Conley indices on the finite dimensional approximations.

\section {Examples}
\label {sec:ex}

Here we apply the results in the previous section to compute explicitely
some Floer spectra. Let $Y$ be a rational homology 3-sphere and $\spc$ a
$spin^c$ structure on it. (In the case when $H_1(Y; \zz) = 0,$ we drop
$\spc$ from the notation, since it is uniquely determined.) Provided we
have a description of the Seiberg-Witten flow lines on $Y,$ we can apply
Proposition~\ref{mhom} to calculate the Floer stable homotopy type
$\swf(Y, \spc).$

The Seiberg-Witten flow is equivariant with respect to the action of the
group $\tt = S^1,$ acting trivially on forms and by rotations on spinors.  
This action is semifree, meaning that there are only free and trivial
orbits. When suspending, we only consider the representations $\rr$ and
$\cc$ of $\tt.$ Correspondingly, there are two type of ``critical
points'':  reducible ones (on which the action is trivial) and
irreducible ones (on which the action is free).

In fact, since $b_1(Y)=0$ there is only one reducible $\theta.$ In 
\cite{M} it is explained how the reducible can be given an absolute Morse 
index equal to $-2n(Y, \spc, g),$ where $n(Y, \spc, g)$ is a combination 
of the Dirac and signature eta invariants for the Riemannian metric which 
we use on $Y.$ If $X$ bounds a 4-manifold with boundary $X,$ and $X$ has 
a $spin^c$ structure with determinant line bundle $\hat L,$ we can express 
$n(Y, \spc, g)$ as:
$$ n(Y, \spc, g) =  \text{ind}_{\cc}(D^+_A) - \frac{c_1(\hat L)^2 -
\sigma(X)}{8} \in \frac{1}{8N} \zz, $$
where $N$ is the cardinality of $H_1(Y; \zz).$ 

Thus, the Conley index of the reducible is a sphere 
$$ S^{-n(Y, \spc, g)\cc} = \bigl( \cc^{-n(Y, 
\spc, g)}\bigr)^+.$$ Note that we allow for negative and even rational 
indices. The choice of an absolute grading for the reducible 
induces one for the irreducibles too. 

Their Conley indices are free cells 
of the form $\Sigma^m (\tt_+).$ Note that we use the notation $\tt$ for a 
circle with free $\tt$ action and $S^1$ for one with trivial $\tt$ 
action; $\tt_+$ stands for $\tt$ together with a disjoint basepoint. We 
should 
also point out that for free cells suspension by $\rr^2$ is the same as 
suspension by $\cc,$ so the notation $\Sigma^m$ is unambiguous.

Finally, note that when we change the orientation on $Y,$ $\swf(- Y, 
\spc)$ is the spectrum dual to $\swf(Y, \spc).$ Morse theoretically, all 
the flow lines go in the reverse direction, and the index of the 
reducible switches sign. For the free cells we can use the 
Wirthm\"uller isomorphism $D(\tt_+) = \Sigma^{-1}(\tt_+),$ which implies 
$D(\Sigma^m(\tt_+)) = \Sigma^{-m-1}(\tt_+).$

\subsection {Elliptic 3-manifolds.}
\label {sec:ell}

The simplest case is that when $Y$ is a quotient of $S^3$ by some finite 
group. Then $Y$ admits a metric of positive scalar curvature $g$, and 
therefore the only critical point is the reducible. The only thing left 
to compute is $n(Y, \spc, g)$ for the metric $g$.

In \cite{M}, we have done this for $S^3$ and for the Poincar\'e sphere 
$P= \Sigma(2,3,5):$
$$ \swf(S^3) = S^0;\hskip6pt \swf(\Sigma(2,3,5)) = \cc^+.$$

When $Y$ is a lens space $L(n,1),$ the invariants $n(Y, \spc, g)$ have 
been computed by Nicolaescu in \cite{N2}. We think of $Y$ as the $S^1$
bundle over $S^2$ of degree $-n$ for $n\geq 1.$ (If we change the
orientation we obtain the bundle of degree $n.$)  Observe that $Y$ bounds
a disk bundle $D(-n) \to S^2.$ The $spin^c$ structures on $D(-n)$ are
denoted $\sps_j, j\in \zz,$ so that $c_1( \det(\sps_j)) = -n +2j \in \zz
\cong H^2(D(-n); \zz).$

Since $H^2(Y; \zz) = \zz/n,$ there are $n$ different $spin^c$ structures
on $Y,$ denoted $\spc_0, \dots, \spc_{n-1},$ such that $\spc_k$ is the
restriction of $\sps_j$ to the boundary for every $j \equiv k \bmod n.$
We have $c_1(\text{det} \spc_k) \equiv 2k \bmod n.$ Then 
$$ \swf(Y, \spc_k) = S^{-n_k \cc},$$
where
$$ n_k = n(Y, \spc_k, g)= \frac{(n-2k)^2 -n}{8n}. $$

\subsection {Some Brieskorn spheres.}

Given $n \geq 1,$ we denote by $Y=\Sigma(2, 3, 6n \pm 1)$ the Brieskorn
sphere oriented as the boundary of the complex singularity $z_1^2 + z_2^3 +
z_3^{6n \pm 1} = 0.$ In this way, for example, $-\Sigma(2, 3, 6n-1)$ can be 
described as the
oriented boundary of the nucleus of the elliptic surface $E(n),$ as well as 
the result of $1/n$ surgery on the right-handed trefoil. On the other hand, 
$-1/n$ surgery on the right-handed trefoil gives $\Sigma(2,3, 6n+1).$

We use the description of the Seiberg-Witten flow lines on $Y$ given by
Mrowka, Ozsv\'ath, and Yu in \cite{MOY}. This description uses a specific
metric and a certain reducible connection on $Y$ instead of the usual
Levi-Civita connection. However, by the continuation properties of the
Conley index, the invariant $\swf$ defined this way is the same as the
usual one. With this choice of metric and connections, all the critical
points are nondegenerate, so the Seiberg-Witten flow is Morse-Bott.
It is usually not Morse-Bott-Smale, and this makes the
application of Proposition~\ref{mhom} more difficult. Nevertheless, we can
still compute the Floer spectra. 

The computation of the absolute index of the reducible in the examples
below was done by Nicolaescu in \cite{N1}. We present the case of $-Y$
rather than $Y$ so that the reader can easily compare the Borel homology of
$\swf$ to the similar computations of $HF^+$ in Ozsv\'ath-Szab\'o theory
(\cite{OS1.5}, \cite{OS2}).

\smallskip

\noindent {\bf \underline{r=12j-1}:} Let us do the case $j=1$ first, which 
will be useful 
to us later. Then the Seiberg-Witten flow on $-Y$ is Morse-Bott-Smale, 
with the reducible of index $0$ and two irreducibles of index $-2.$ Modulo 
the action of $\tt,$ there
is exactly one flow line from the reducible to each irreducible.
It follows that there is an exact triangle in stable homotopy:
\begin {equation}
\label {ya}
\cdots \to S^{-1} \to \Sigma^{-2}(\tt_+) \vee \Sigma^{-2}(\tt_+) \to  
\swf(- Y) \to S^0 \to \cdots
\end {equation}

The connecting morphism, corresponding to the attaching map from the
trivial 2-cell to the free cells, can be thought of as an element in a
stable cohomotopy group:
$$ \{ S^{-1}, \Sigma^{-2}(\tt_+ \vee \tt_+) \}^{\tt} =
\{\Sigma^2 D(\tt_+
\vee \tt_+), S^1\}^{\tt} = $$
$$ =\pi^0_{\tt}(\tt_+ \vee \tt_+) = \pi^0(S^0) \oplus \pi^0(S^0) = \zz
\oplus \zz.$$

Given the count of flow lines, this morphism must be $(\pm 1, \pm 1).$
The signs depend on the (noncanonical) identification of the Conley 
indices of the critical points with the standard cells $S^0, 
\Sigma^{-2}\tt_+.$ It is not important to fix a convention at this point. 

Nonequivariantly, we can identify $\tt_+$ with $S^1 \vee S^0$ by
choosing a basepoint on $\tt.$ The triangle (\ref{ya})
becomes:
$$ \cdots \to S^{-1} \to S^{-1} \vee S^{-2} \vee S^{-1} \vee S^{-2} \to 
\swf(-Y) \to S^0 \to \cdots. $$
This shows that, nonequivariantly, there is an isomorphism:
$$ \swf(-\Sigma(2,3,11)) = S^{-2} \vee S^{-2} \vee S^{-1}.$$

More generally, for $j > 1,$ there is one reducible of index $0$ and $2j$
irreducibles of index $-2.$ Again, there is one flow line from the
reducible to each irreducible. However, there are also flow lines between
different irreducibles, so the flow is not Morse-Bott-Smale. It turns out
that we do not need to understand the topology of the spaces of flow lines
in this case. Indeed, we can form attractor-repeller pairs by splitting
along the levels of the CSD functional. Each irreducible contributes 
$\Sigma^{-2}(\tt_+),$ and whenever we attach two of them together we do so 
by a morphism:
$$ \Sigma^{-2}(\tt_+) \to \Sigma^{-1}(\tt_+).$$

All these morphisms must be trivial, which proves that all the irreducibles 
contribute with a wedge of the respective cells. It follows that for 
$Y=\Sigma(2,3, 12j-1),$ the Floer spectrum admits a presentation:
$$ \cdots \to S^{-1} \to \bigvee_{i=1}^{2j} \Sigma^{-2}(\tt_+)  \to
\swf(- Y) \to S^0 \to \cdots $$
with the connecting morphism being $(\pm 1, \dots, \pm 1) \in \zz^{2j}.$

\smallskip
\noindent {\bf \underline{r=12j-5}:} This case is similar to the
previous one, except for a shift in index: there is one reducible of 
complex index $2$ and $2j$ irreducibles of index $0.$ We get:
$$ \swf(-\Sigma(2,3,12j-5)) = \Sigma^{\cc} \swf(-\Sigma(2,3,12j-1)).$$

\smallskip
\noindent {\underline{\bf r=12j+1}:} There is one reducible and $2j$ 
irreducibles, all 
of index $0.$ The flow is not Morse-Smale-Bott even for $j=1,$ because 
there are flow lines from the reducible to the irreducibles. The 
Chern-Simons-Dirac functional $CSD$ 
takes a bigger value on the reducible than on any irreducible, so we 
can find an attractor-repeller exact triangle by splitting along a level 
set of $CSD:$
$$ \cdots \to S^{-1} \to \bigvee_{i=1}^{2j} \tt_+ \to \swf(- Y) \to S^0 
\to \cdots $$

The connecting morphism must be zero, simply because:
$$ \pi_{-1}^{\tt} (\tt_+) = \pi^2_{\tt}(\tt_+) = 
\pi^2(S^0) =0.$$

Therefore, we have even equivariantly:
$$\swf(-\Sigma(2,3,12j+1)) = S^0 \vee \bigl(\bigvee_{i=1}^{2j} 
\tt_+\bigr).$$

\smallskip
\noindent {\bf \underline{r=12j+5}:} This case is similar to the previous 
one, except 
all the critical points have index $-2.$ We get:
$$ \swf(-\Sigma(2,3,12j+5)) = \Sigma^{-\cc} \swf(-\Sigma(2,3,12j+1)).$$

\section {Applications}

\subsection{Adjunction inequalities.}

We start with some applications of Theorem~\ref{gluing} to gluings along
lens spaces. The two propositions proved in this section are not new (they
appear in \cite{FS}), but it is interesting to see how they can be
obtained with our techniques.

Consider now a closed, orientable 4-manifold $X$ with $\pi_1(X)=1.$ We can
identify the set of $spin^c$ structures $\sps$ on $X$ with the set of
characteristic elements $c \in H^2(X; \zz)$ via the correspondence $\sps
\to c=c_1(\hat L).$

Recall that an element $c \in H^2(X; \zz)$ 
is called a {\bf basic class} if the Seiberg-Witten invariant $\sw(X, 
\sps) \neq 0.$ The manifold $X$ is called {\bf of simple type} if all 
basic calsses satisfy $c^2 = 3\sigma(X) +2\chi(X),$ or equivalently, 
if $\sw=0$ whenever the formal dimension of the monopole moduli space is 
nonzero. 

It turns out that the presence of embedded surfaces in specific homology
classes of our 4-manifold $X$ imposes some restrictions on the set of
basic classes. These restrictions come from the so-called ``adjunction
formulae.'' The most general version of these formulae appears in
\cite{OS1}. Our tools are insufficient for proving these results in the
general case, but they allow us to find the constraints imposed by the
existence of embedded \emph {spheres.}

Let $\Sigma \subset X$ be an embedded sphere. We study the case
$[\Sigma]^2 = N > 0$ first. In this situation a neighborhood of $\Sigma$
is the disc bundle $D(N)$ over $S^2$ with boundary the lens space $L(N,
N-1),$ which is $L(N,1)$ with the opposite orientation. This gives a
decomposition of $X$ into two pieces $X'$ and $D(N),$
glued along their common boundary $L(N,N-1).$ We have $b_2^+(D(N)) =1$ and
$\swf(L(N, 1), \spc) \cong S^{-n\cc},$ where $n$ depends on $\spc.$ 
Hence
the relative invariant $\Psi(D(N))$ lives in 
$ \pi^1_{\tt}(S^{d\cc})$ for some $d\in \zz.$ These groups are 
torsion for every $d.$ Since 
$\Psi(D(N))$ is torsion, so  
is $\Psi(X)$ by virtue of Theorem~\ref{gluing}, so the Seiberg-Witten 
invariant is zero. We conclude:  

\begin {proposition}
Let $X$ be a smooth, closed, oriented, simply connected 4-manifold with
$b_2^+(X) >1.$ If there exists an embedded sphere $\Sigma \subset X$ with
$[\Sigma]^2 >0,$ then $X$ has no basic classes.
\end {proposition}

Let us now consider the case $[\Sigma]^2 = -N <0,$ under the additional
assumption that $X$ is of simple type. Let $c$ be a basic class of $X,$   
with corresponding $spin^c$ structure $\sps.$ A neighborhood of $\Sigma$
is the disc bundle $D(-N)$ over $S^2$ with boundary the lens space $L(N, 
1).$ This decomposes
$X$ into two pieces $X'$ and $D(-N)$ and breaks $\sps$ into $spin^c$
structures $\sps'$ and $\sps_j$ on $X'$ and $D(-N),$ respectively. 
We must have $-N+2j = c([\Sigma]).$ With the notations from 
Subsection~\ref{sec:ell}, the induced
$spin^c$ structure on $L(N, 1)$ is $\spc_k,$ with $0 \leq k < N-1,
k\equiv j \bmod N.$ Note that $b_2^+(D(-N))=0,$
so the invariant $\Psi(D(-N), \sps_j)$ lives in $\pi^0_{\tt} 
(S^{i\cc}),$ where
$$i=\ind_{\cc} (D_{\hat A}^+) = n_k +\frac{-(2j-N)^2/N+1}{8}.$$

Recall that $n_k= ((N-2k)^2-N)/8N.$ It follows that:
$$ i= \frac{(N-2k)^2 - (N-2j)^2}{8N}.$$

The following lemma appears in \cite{B}:
\begin {lemma}
\label {d2}
Let $f: (\rr^m \oplus \cc^{n+d})^+ \to (\rr^m \oplus \cc^n)^+$ be a
$\tt$-equivariant map such that the induced map on the fixed point sets
has degree $1.$ Then $d \leq 0$ and $f$ is $\tt$-homotopic to the
inclusion.
\end {lemma}

Therefore, $\Psi(D(-N), \sps_j)$ is the class of the
inclusion. We claim that $i=0.$ Indeed, if $i \neq
0$ and $\sw(X, \sps) \neq 0,$ we could consider the $spin^c$ structure
$\sps_{new}$ on $X$ obtained from $\sps'$ on $X'$ and $\sps_k$ (instead of
$\sps_j$) on $D(-N).$ By two applications of the gluing theorem we would
get $\sw(X, \sps_{new}) =\sw(X, \sps) \neq 0.$ But $c(\det(\sps_{new}))^2
= c(\det(\sps))^2 -i,$ so $X$ would not be of simple type.

Therefore, we must have $i=0,$ or $(N-2k) = \pm (N-2j).$ This happens if   
and only if $|N-2j| \leq N.$ We deduce the following:

\begin {proposition}
Let $X$ be a smooth, closed, oriented, simply connected 4-manifold with
$b_2^+(X) >1.$ If $X$ is of simple type and there exists an embedded
sphere $\Sigma \subset X$ with $[\Sigma]^2 = -N < 0,$ then every basic
class $c$ of $X$ satisfies $|c([\Sigma])| \leq N.$
\end {proposition}

\subsection {Exotic nuclei.} Here we present the proof of 
Theorem~\ref{exotic}. The Bauer-Furuta invariants of $X = K3 \# K3 
\# K3$ have been computed in \cite{B}. They are zero except 
for the trivial $spin^c$ structure $\spc_0,$ when, nonequivariantly:
$$ \Psi(X, \spc_0) = 12 \in \pi_3(S^0) = \zz/24.$$

Set $Y= -\Sigma(2,3,11),$ oriented now as the boundary of the 
nucleus $N(2) \subset E(2)=K3.$ Assume that $X$ decomposes as $X_1 
\cup_Y X_2,$ where $X_1$ is an exotic nucleus $N(2)_{p,q}$ with $(p,q) 
\neq (1,1),$ and $X_2$ has intersection form $6(-E_8) \oplus 8H.$ 

Recall from Section~\ref{sec:ex} that, nonequivariantly:
$$ \swf(-Y) = S^2 \vee S^2 \vee S^1.$$ 

Let $\spc'$ be the restriction of $\spc_0$ to $X_2.$ Then the 
nonequivariant Bauer-Furuta invariant of $(X_2, \spc')$ lives in:
$$ \pi_4(\swf(-Y)) = \zz/2 \oplus \zz/2 \oplus \zz/24.$$
Denote by $z$ its projection to the $\zz/24$ factor.

Dually, we have that $\swf(Y) = S^{-2} \vee S^{-2} \vee S^{-1}$ 
nonequivariantly. The Bauer-Furuta invariant of $X_1$ lives in:
$$\pi_{-1}(\swf(Y)) = \zz/2 \oplus \zz/2 \oplus \zz.$$

We can get additional information by observing that it has to lie in the 
image of the equivariant group $G=\pi^{\tt}_{-1}(\swf(Y))$ under the 
natural forgetting map. To compute $G$, we look at the triangle 
($\ref{ya}$):
$$\cdots \to S^{-1} \to \Sigma^{-2}(\tt_+) \vee \Sigma^{-2}(\tt_+) \to 
\swf(Y) \to S^0 \to \cdots $$

Note that:
$$ \pi_{-1}^{\tt}(\Sigma^{-2} (\tt_+)) = \pi^1_{\tt}(D(\Sigma^{-2} 
(\tt_+))) = \pi^0_{\tt}(\tt_+) = \zz,$$
and the forgetting map to $\pi_{-1}(\Sigma^{-2} (\tt_+)) = \zz/2 \oplus 
\zz$ is given by $1 \to (0,1).$

Applying the functor $ \pi_{-1}^{\tt}$ to (\ref{ya}) we get a long exact 
sequence:
$$\cdots\to \zz \to \zz \oplus \zz \to \pi_{-1}^{\tt}(\swf(Y)) \to 0$$
with the first map being the diagonal $1 \to (1,1).$ 

Consequently, we find that $G = \pi_{-1}^{\tt}(\swf(Y)) = \zz$ and the 
forgetting map to $\zz/2 \oplus \zz/2 \oplus \zz$ is given by $1 \to 
(0,0,1).$

Let $P$ be the primitive class in $H_2(N(2)_{p,q}; \zz)$ such
that $pqP = F$ is the class Poincar\'e dual to the 
elliptic fiber. For a
$spin^c$ structure $\spc=mP, m\in \zz$ on $X_1 = N(2)_{p,q},$ the 
Bauer-Furuta invariant 
$$\Psi(X_1, \spc) \in G= \zz$$ 
is equal to the relative Seiberg-Witten invariant $SW(X_1, \spc)$ as 
defined in \cite{SS}. This counts the difference in the number of monopoles on 
$X_1$ which restrict on the boundary to each of the two reducibles. The formal 
relative Seiberg-Witten series was computed by Stipsicz and Szab\'o to be:
$$ SW(X_1) = \frac{\sinh^2(pqP)}{\sinh(pP) \cdot 
\sinh(qP)}.$$

Let
$$x_0 = SW(X_1, 0); \hskip6pt x_1=SW(X_1, (2pq-p-q)P)=1.$$

Build the $spin^c$ structure $\spc_1$ on $X$ by gluing $\spc'$ 
to $(2pq-p-q)P.$ Nonequivariantly, Theorem~\ref{gluing} says  the 
Bauer-Furuta invariant of $X$ with either $\spc_0$ and $\spc_1$ is 
obtained from those of $X_1$ and $X_2$ via the duality map:
\begin {eqnarray*}
\pi_{-1}(\swf(Y)) \times  \pi_4(\swf(-Y))  &\to& \pi_3(S^0) \\
(\zz/2 \oplus \zz/2 \oplus \zz) \times (\zz/2 \oplus \zz/2 \oplus \zz/24) 
&\to& \zz/24
\end {eqnarray*}
given by:
$$ (a_1, b_1, c_1) \times (a_2, b_2, c_2) \to \bigr( 12(a_1 a_1' + b_1 
b_1') + c_1 c_2\bigr) \text{ mod }24.$$

Therefore:
$$\Psi(X, \spc_0) = x_0z \in \zz/24,$$
while
$$ \Psi(X, \spc_1) = x_1z =z  \in \zz/24.$$ 

On the other hand, we knew from the beginning that $\Psi(X, \spc_0) = 12,$ 
while $\Psi(X, \spc_1) = 0.$ This is a contradiction. Hence $X$ does not 
contain an exotic nucleus.

\end{document}